\documentclass[]{interact}

\usepackage{amsthm}
\usepackage{amssymb} 
\usepackage{amsmath}
\usepackage{fancyhdr}
\usepackage{hyperref}
\usepackage{mathtools}
\usepackage{physics}
\usepackage{color}

\numberwithin{equation}{section}

\newtheorem{definition}{Definition}[section]
\newtheorem{theorem}[definition]{Theorem}
\newtheorem{lemma}[definition]{Lemma}
\newtheorem{proposition}[definition]{Proposition}
\newtheorem{corollary}[definition]{Corollary}

\title{On the topological entropy of saturated sets for amenable group actions}

\begin{document}

\bibliographystyle{plain}

\author{
	\name{Xiankun Ren\textsuperscript{a,$^\ast$} Xueting Tian\textsuperscript{b,$^{\S}$} and Yunhua Zhou\textsuperscript{a,$^{\ddag}$}
		\thanks{\textsuperscript{$^\ast$}Email:xkren@cqu.edu.cn; \textsuperscript{$^{\S}$}xuetingtian@fudan.edu.cn; \textsuperscript{$^{\ddag}$}zhouyh@cqu.edu.cn
		\newline \indent Corresponding author: Xueting Tian
		\newline \indent2010 Mathematics Subject Classification.  37B05,37B40,54H20
}
	}	
	\affil{\textsuperscript{a}College of Mathematics and Statistics\\Chongqing University
		\\ Chongqing 401331, China}	
	\affil{\textsuperscript{b}School of Mathematical Sciences, Fudan University
		\\ Shanghai 200433, China}
}

\maketitle

\begin{abstract}
Let $(X,G)$ be a $G$-action topological system, where $G$ is a countable infinite discrete amenable group and $X$ a compact metric space. We prove a variational principle for topological entropy of saturated sets for systems which have the specification property and uniform separation property. We show that certain algebraic actions satisfy these two conditions. We give an application in  multifractal analysis.
\end{abstract}


\section{Introduction}

\indent In this paper, a dynamical system $(X, G)$  always means that $X$ is a compact metric space and $G$ a countable discrete amenable group acting on $X$ continuously. Let $M(X)$ stand for  the set of all Borel probability measures endowed with weak$^{*}$ topology, $M(X,G)\subset M(X)$  stand for the set of $G$-invariant measures and $E(X,G)\subset M(X)$ be the set of ergodic $G$-invariant measures.

We are interested in comparing the metric entropy of $\mu\in M(X,G)$ with the topological entropy, which is a measure of complexity of the dynamical system. In this paper, we are dealing with general group actions instead of $\mathbb{Z}$-actions. The problem is interesting because new phenomena and difficulties arise as we go to more general group actions.   In 1987, Ornstein and Weiss \cite{ornstein1987entropy} developed the so-called {\itshape quasi-tiling} method, which has been a basic tool in the study of amenable group actions. The quasi-tiling  can serve as the substitute of the Rokhlin tower and allows people to generalize the known results for $\mathbb{Z}$-actions to amenable group actions. Many people have made lots of progress in many directions of group actions. For example, Kieffer \cite{kieffer1975generalized} extended the definition of metric entropy to a probability measure preserving amenable group action and Ornstein and Weiss \cite{ornstein1987entropy}  generalized the Ornstein theory to this setting. The Shannon-McMillan theorem for amenable group actions is due to Kieffer  \cite{kieffer1975generalized}. In the ergodic case, Ornstein and Weiss obtained almost everywhere convergence for a special type of F{\o}lner sequence that exits for some amenable groups but not others \cite{ornstein1983shannon}. Lindenstrauss \cite{lindenstrauss1999pointwise} later established this pointwise result for tempered F{\o}lner sequence with superlogarithmic growth, which can be found in any amenable group. The variational principle for amenable group actions was due to Ollagnier \cite{ollagnier2007ergodic}. The pointwise ergodic theorem for amenable group actions is due to Lindenstrauss \cite{lindenstrauss1999pointwise}. But there are still many important theory in integer actions are not confirmed for general group actions. In this paper we try to confirm some such kind theory in amenable group actions. Expansiveness and specification property are needed in our setting. Chung and Li \cite{chung2015homoclinic} extended the definition of specification property to general group actions. There are non-trivial examples under general actions with expansiveness and specification, see \cite{ren2018}. For more information of amenable group actions, readers may refer \cite{huang2011local,ollagnier2007ergodic,kerr2016ergodic}.

The study of the thermodynamic formalism and multifractal analysis for maps with some hyperbolicity has drawn the attention
of many researchers from the theoretical physics and mathematics communities in the last decades. The general concept of
multifractal analysis, which can be traced back to Besicovitch, is to decompose the phase space  in subsets of points which have a similar dynamical behavior and to describe the size of each of such subsets from the geometrical or topological viewpoint.
We refer the reader to \cite{olsen1995multifractal,pesin1997multifractal} and lots of such kind progress \cite{barreira2000sets,barreira2008dimension,takens2003variational,olsen1995multifractal,olsen2003normal, barreira1997multifractal, todd2008multifractal, climenhaga2010multifractal, gelfert2007lyapunov, liang2017variational, jordan2008multifractal, bomfim2017multifractal, thompson2012irregular,varandas2012non,barreira2000multifractal} etc. and references therein under different settings.  It is still necessary to let people know which multifractal analysis of general group actions hold.  It is almost not possible for us to generalize so many results in multifractal analysis for $\mathbb{Z}$ actions one by one to general group actions since this is  a huge hard work. Observe that the study of saturated sets is the most critical  technique  which can imply various results in multifractal analysis including irregular sets, level sets and classification of recurrent and transitive points by constructing different saturated sets of measures with some required information, for example, see \cite{pfister2007topological,tian2016different,huang2019transitively}. Thus in this paper we will give such a charaterization for topological entropy on saturated sets so that one can use this result to  know which kind multifractal analysis to be get for  amenable group actions. Here we can not give all the applications but give some applications,  for example, we will give some application on irregular sets and level sets of continuous observables. This result may have applications in various multifractal analysis, including classification of transitive points\cite{tian2016different,huang2019transitively}, level sets and irregular sets of  asymptotically additive or almost additive continuous observables\cite{feng2010lyapunov,barreira2009almost} and their mixed version \cite{barreira2001variational} or higher version \cite{barreira2002higher}, etc.

The set of invariant measures plays an important role in the study of ergodic theory. For $\mathbb{Z}$-actions, the invariant measure always exists. But for general group actions $(G,X)$, the set of $G$-invariant measures may be empty. A well known result shows that when $G$ is amenable, there always exists a $G$-
invariant measure. The class of amenable group includes all finite groups, Abelian groups and solvable groups.

Let $\mathcal{F} = \{F_n\}$ be a F{\o}lner sequence. We will study the empirical measure (along $F_n$) of $x$, which is the probability measure defined as
\[
\mathcal{E}_{F_n}(x) := \frac{1}{|F_n|}\sum_{s\in F_n}\delta_{sx},
\]
where $\delta_x$ is the Dirac mass at $x\in X.$

 A subset $D\subset X$ is called {\itshape saturated} with respect to $\mathcal{F}$ if $x\in D$ and the sequences $\{\mathcal{E}_{F_n}(x)\}$
and $\{\mathcal{E}_{F_n}(y)\}$ have the same limit-point set, then $y\in D.$ The limit point set of  $\{\mathcal{E}_{F_n}(x)\}$ is always a non-empty compact subset $V(x,\mathcal{F})\subset M(X).$  For each non-empty closed subset $K\subset M(X)$, denote by $G_K(\mathcal{F}) := \{x\in X\mid V(x,\mathcal{F}) = K\}$ the {\itshape saturated set} of $K$ with respect to $\mathcal{F}$. 

The existence of saturated sets was firstly showed by Sigmund for systems with specification including hyperbolic systems \cite{sigmund1972space,sigmund1974dynamical} and then was generalized to non-uniformly hyperbolic and non-uniformly expanding systems  \cite{liang2010ergodic,tian2017topological}. This result can imply that the points whose empirical measures equal to the space of invariant measures form a residual set and in particular, every irregular set is either empty or residual in the whole space, see \cite{tian2016different}.

The entropy estimate on saturated sets was firstly studied in \cite{pfister2007topological} for systems with specification including hyperbolic systems and then was generalized  recently to non-uniformly hyperbolic and non-uniformly expanding systems \cite{liang2017variational,tian2017topological}. In particular, remark that the entropy estimate of the particular saturated set of ergodic measures is  due to Bowen \cite{bowen1973topological}. Here we give the existence and entropy estimate of saturated sets for amenable group actions. For $\mathbb{Z}$-actions, a saturated set is always invariant but for amenable group actions, a saturated set may not be $G$-variant. We remark that we do not need saturated sets to be $G$-invariant in Theorem \ref{mainthm1}.

\begin{theorem}\label{mainthm1}
	Let $(X,G)$ be dynamical system satisfying
	 the specification property and uniform separation property. Let $\mathcal{F}=\{F_n\}$ be a F{\o}lner sequence with $\frac{|F_n|}{\log n} \rightarrow \infty$ then for any non-empty connected closed subset $K\subset M(X,G)$
	\[
	h^{B}_{top}(G_K(\mathcal{F}),\mathcal{F}) = \inf \{h_{\mu}(X,G)\mid \mu\in K\}.
	\]
	Here $h^{B}_{top}(G_K(\mathcal{F}),\mathcal{F})$ is the Bowen topological entropy of $G_K(\mathcal{F})$ with respect to $\mathcal{F}$.
\end{theorem}

We point out that the set $V(x,\mathcal{F})$ may not be in general connected (it depends on the properties of $\mathcal{F}$, see \cite[Page 694]{lkacka2018quasi}). But we still need the set $K$ to be connected. The reason is that we need two levels of standard bricks to cover $F_n$ (see Lemma \ref{section3lemma2}) and as a result we need the sequence of measures $\{\alpha_n\}$ in Section 5.1 satisfying $\lim\limits_{n\rightarrow\infty}D(\alpha_n,\alpha_{n+1}) =0.$

{\itshape Remark 1.1:} In \cite{zhang2018topological}, the author introduced $g-$almost product property for amenable group actions and showed that specification property implies $g-$almost product property. Our proof of Theorem \ref{mainthm1} can be modified to $g-$almost product property cases by reconstructing $Y_K$ using different separated points with respect to $g-$almost product property. There is no examples to show the difference between specification property and $g-$almost product property for group actions yet. Here for simplicity of the writing, we just prove the cases which satisfy the specification property.

{\itshape Remark 1.2:} For $K$ being a singleton, Theorem \ref{mainthm1} means the existence of generic points. This would extend a result of Lacka \cite{lkacka2017generic} who (extending an old result of Sigmund \cite{sigmund1977on} and completing a result
of Ren \cite{ren2018}) proved that every invariant measure  for an amenable residually finite group action satisfying the weak
specification property has a generic point.
\subsection{Applications to multifractal analysis }

From the measure theoretical viewpoint, the pointwise ergodic theorem guarantees that with respect to a tempered F{\o}lner sequence the irregular set has zero measure for every invariant measure. Nevertheless, irregular sets may have full topological entropy, see \cite{barreira2000sets,thompson2012irregular,dong2015irregular,liang2017variational}. The classical approach to prove that the irregular set of continuous observables for   $\mathbb{Z}$ action that are not cohomologous to a constant has
full topological entropy uses the uniqueness of equilibrium states\cite{barreira2001variational,barreira2002higher,barreira2000sets}. However, up to now it still unknown the uniqueness of equilibrium states for amenable group actions. Here we construct different saturated sets to get the role. However, two problems are still unknown that whether the spectrum of level sets has some smoothness with respect to the level and whether there is an ergodic measure supported on the level set with metric entropy same as the topological entropy of the level set.

Let $\varphi \in C(X,\mathbb{R})$ and $\mathcal{F}=\{F_n\}$ be a F{\o}lner sequence, then $X$ can be divided into the following parts:
\[ X = \bigcup_{\alpha \in \mathbb{R}}X(\varphi, \alpha,\mathcal{F}) \cup \hat{X}(\varphi,\mathcal{F})\]
where for $\alpha\in \mathbb{R},$
$$X(\varphi, \alpha,\mathcal{F}) = \{x\in X \mid \lim\limits_{n\rightarrow \infty}\frac{1}{|F_n|}\sum\limits_{s\in F_n}\varphi(sx) = \alpha\}$$
and
$$\hat{X}(\varphi,\mathcal{F}) = \{x\in X \mid \lim\limits_{n\rightarrow \infty}\frac{1}{|F_n|}\sum\limits_{s\in F_n}\varphi(sx) \text{ does not exist} \}.$$

The set $X(\varphi, \alpha,\mathcal{F})$ is called a {\itshape level set} with respect to $\mathcal{F}$ and $\varphi$ and the set $\hat{X}(\varphi,\mathcal{F})$ is called {\itshape the historic set} with respect to $\mathcal{F}$ and $\varphi$ or $\varphi$-{\itshape irregular set.} 

The level set is called the {\itshape multifractal decomposition set} of ergodic average of $\varphi.$ In particular, one is interested in the ‘size’ of these sets $X(\varphi, \alpha,\mathcal{F}).$ For irregular sets, Pesin and Pitskel \cite{pesin1984topological} are the first to notice the phenomenon of the irregular set carries full topological entropy in the case of the full shift on two symbols from the dimension perceptive. Ruelle \cite{rulle2001} used the
terminology "historic behavior" to describe irregular points and in contrast to dimensional perspective.

For amenable group actions, we have the following results.
\begin{theorem}\label{mainthm2}
Let $(X,G)$ be a dynamical system and $\mathcal{F}=\{F_n\}$ be a F{\o}lner sequence with $\frac{|F_n|}{\log n}\rightarrow \infty$. Suppose the system has the specification and uniform separation properties. If $\varphi\in C(X,\mathbb{R})$ and $\hat{X}(\varphi,\mathcal{F})$ is non-empty, then
		\begin{align}
		h_{top}^B(\hat{X}(\varphi,\mathcal{F}),\mathcal{F}) = h_{top}(X,G).
			\end{align}
\end{theorem}

\begin{theorem}\label{mainthm3}
	Let $(X,G)$ be a dynamical system and $\mathcal{F}=\{F_n\}$ be a F{\o}lner sequence with $\frac{|F_n|}{\log n}\rightarrow \infty$. Suppose the system has the specification and uniform separation properties.
For $\alpha\in \mathbb{R}$ and $\varphi\in C(X,\mathbb{R})$ we have
			\begin{align}
		h^B_{top}(X(\varphi, \alpha,\mathcal{F}),\mathcal{F}) = \sup\{h_{\mu}(X,G)\mid \mu\in M(X,G), \int_{X}\varphi \dd\mu = \alpha\}.
			\end{align}
\end{theorem}

This paper is organized as follows. In Section 2, we recall some definitions and in Section 3 we define uniform separation property for amenable group actions. In Section 4, we prove the upper bound of Theorem 1.1 and in Section 5 we prove the lower bound of Theorem 1.1. In Section 6, we prove Theorem \ref{mainthm2} and Theorem \ref{mainthm3}. In the Appendix, we prove Theorem \ref{approthm}.

\section{Preliminaries}

In this section, we will introduce some notions and properties.

\subsection{Metric on $X$ and $M(X)$}
Let $\varphi\in C(X,\mathbb{R})$ and $\mu\in M(X)$. We set 
\[
\langle\varphi,\mu\rangle = \int_X \varphi  \dd \mu.
\]
There exists a countable  separating set of continuous functions $\{\varphi_1,\varphi_2,\dots\}$ with $0\leq \varphi_k\leq 1,$ and such that 
\[
D(\mu,\nu) := \sum_{k=1}^{\infty}2^{-k}|\langle \varphi_k,\mu\rangle - \langle \varphi_k,\nu\rangle|,
\]
defines a compatible metric for the weak$^*$-topology on $M(X).$ For $r>0$ and $\mu\in M(X)$, define
\[
\mathcal{B}(\mu,r) = \{\nu\in M(X)\mid D(\mu,\nu) < r\}.
\]
In this paper, for convenience we will use an equivalent metric 
\begin{align}
\rho(x,y) := D(\delta_x,\delta_y). \label{metriconx}
\end{align}
as the metric on $X$.

\subsection{Amenable groups and tilings of amenable groups}

Let $F(G)$ be the collection of finite subsets of $G$. 

\begin{definition}
	Let $\Omega, K\in F(G)$ be two finite subsets of a group $G.$ The $K$-interior of $\Omega$ is the subset $Int_{K}(\Omega)$ defined by 
	\[
	Int_{K}(\Omega) := \{g\in G \mid Kg \subset \Omega\}.
	\]
The $K$-closure of $\Omega$ is the subset $Cl_{K}(\Omega) \subset G$ defined by
	\[
	Cl_{K}(\Omega) := \{g\in G\mid Kg \cap \Omega \neq \emptyset\}.
	\]
The $K$-boundary of $\Omega$ is the subset $\partial_{K}(\Omega)\subset G$ defined by
	\[
	\partial_{K}(\Omega) := Cl_{K}(\Omega)\setminus Int_{K}(\Omega).
	\]
The relative amenablity constant of $\Omega$ with respect to $K$ is the number $\alpha(\Omega,K)$ defined by
	\[
	\alpha(\Omega,K) := \frac{|\partial_{K}(\Omega)|}{|\Omega|}.
	\] 
We say	$\Omega$ is  $(K,\delta)$-invariant if $\alpha(\Omega,K) < \delta.$
\end{definition}

A sequence  $\{F_n\}\subset F(G)$ is called a {\itshape F{\o}lner} sequence if for any $s\in G,$
\[
\lim\limits_{n\rightarrow \infty}\frac{|sF_n\triangle F_n|}{|F_n|} = 0.
\] 
We say $G$ is {\itshape amenable}, if it admits a F{\o}lner sequence. Note that if $\mathcal{F}=\{F_n\}$ is a F{\o}lner sequence, then for every $\varepsilon>0$ and $K\in F(G)$ there is $N\in\mathbb{N}$ such that $F_n$ is $(K,\varepsilon)$-invariant for every $n\geq N.$ By $\mathcal{F}$ we always denote a F{\o}lner sequence $\{F_n\}$.

A F{\o}lner sequence $\mathcal{F}=\{F_n\}$ is {\itshape tempered} if for some $C>0$ and all $n\in\mathbb{N}$ one has $|\bigcup_{k\leq n}F^{-1}_{k}F_{n+1}|< C|F_{n+1}|$.

The quasi-tiling-theory is a useful tool for amenable group actions which is set up by Ornstein and Weiss in \cite{ornstein1987entropy}.

Subsets $A_1, A_2,\cdots, A_k\in F(G)$ are {\itshape $\var$-disjoint} if there exists $\{B_1,B_2,\cdots,B_k\}\subset F(G)$ such that

\begin{itemize}
	\item $ B_i\subset A_i \quad i=1,2,\cdots,k,$
	\item $B_i\cap B_j = \emptyset \quad 1\leq i< j\leq k,$
	\item  $\frac{|B_i|}{|A_i|}>1-\var\ \quad i=1,2,\cdots,k.$
\end{itemize}

\noindent For $\alpha\in (0,1],$ we say $\{A_1,A_2,\cdots,A_k\}$ {\itshape $\alpha$-covers} $A\in F(G)$ if
$$\frac{|A\cap(\bigcup_{i=1}^{k}A_i)|}{|A|}\geq \alpha.$$

\noindent We say that $\{A_1,A_2,\cdots,A_k\}\subset F(G)$ is a \ {\itshape $\var$-quasi-tile} of $A\in F(G)$ if there exists $\{C_1,C_2,\cdots,C_k\}\subset F(G)$ satisfying

\begin{itemize}
	\item $A_iC_i\subset A $ and $\{A_ic\mid c\in C_i\}$ forms a $\delta$-disjoint family for $i=1,2,\cdots,k,$
	\item $A_iC_i\cap A_jC_j = \emptyset \quad 1\leq i\neq j\leq k,$
	\item $\{A_iC_i\mid i=1,2,\cdots k\}$ forms a $(1-\delta)-$cover of $A.$
\end{itemize}
Such $C_1,C_2,\cdots,C_k$ are called the {\itshape tiling centers}.

The following proposition is a fundamental quasi-tiling property of amenable groups. The description is a little bit different from \cite[Theorem 6 in I.2]{ornstein1987entropy}, but the ideas are the same.

\begin{proposition}{\cite[Lemma 9.4.14]{coornaert2015topological}}\label{oldtiling}
	Let $G$ be a group and $0< \varepsilon \leq \frac{1}{2}$. Then there exists an integer $s_0 = s_0(\varepsilon) \geq 1$ such that for each $s \geq s_0$ the following holds.
	If $K_1, K_2,\dots, K_s$ are non-empty finite subsets of $G$ such that 
	\[
	\alpha(K_k,K_j)\leq \varepsilon^{2s} \text{ for all } 1\leq j< k \leq s,
	\]
	and $D$ is a non-empty finite subset of $G$ such that
\[
\alpha(D,K_j) \leq \varepsilon^{2s} \text{ for all } 1\leq j \leq s,
\]
then $D$ can be $\varepsilon-$quasi tiled by $K_1, K_2, \dots, K_s.$
\end{proposition}

{\bfseries Remark 2.1:} Let $K_1, K_2, \dots, K_s$ be an $\varepsilon-$quasi-tile of $D \subset G$ and $\{C_j, j=1,2,\dots,s\}$ be the tiling centers. We can modify the tile to get a $(1-\varepsilon)^2-$disjoint cover of $D$ by shrinking every translation of $K_i, \ i = 1, 2, \dots, s.$ In fact, for each $j$, since $\{K_j c_j \mid c_j\in C_j\}$ are $\varepsilon-$disjoint, we can choose $K_j(c_j) \subset K_j$ with $\frac{|K_j(c_j)|}{|K_j|}\geq 1-\varepsilon$ and the elements in $\{K_j(c_j)c_j\}$ are pairwise disjoint. Thus elements in the collection $\{K_j(c_j)c_j\mid c_j\in C_j, \ j=1,2,\dots, s\}$ are pairwise disjoint and 
\[
\frac{|\cup_{j=1}^{s}\cup_{c\in C_j}K_j(c_j)c_j|}{|D|} \geq (1-\varepsilon) \frac{|\cup_{j=1}^{s}\cup_{c_j\in C_j}K_jc|}{|D|} \geq (1-\varepsilon)^{2}.
\]

Also we need a tiling result in \cite{downarowicz2019tilings}.

\begin{definition}
	We say $\mathcal{T}$ is a {\itshape tiling} of $G$ if there exist a {\itshape shape} set $\mathcal{S} =\{S_i \in F(G) \mid 1\leq j \leq k\}$ and tiling centers $C_1, C_2,\dots, C_k$ such that 
	\[
	\mathcal{T}:=\{S_jg \mid g\in C_j, j=1,2,\dots,k\}
	\]
	with $G = \cup\mathcal{T}$ and $A\cap B =\emptyset$ for $A\neq B \in\mathcal{T}.$ Let $\{\mathcal{T}_k\}_{k\geq 1}$ be a sequence of tilings of $G$, we say $\{\mathcal{T}_k\}_{k\geq 1}$ is {\itshape congruent} if for each $k\geq 1$, each element in $\mathcal{T}_{k+1}$ is a union of elements in $\mathcal{T}_k.$
\end{definition}

The following lemma is part of \cite[Lemma 5.1]{downarowicz2019tilings}.

\begin{lemma}\label{newtiling}
	Fix a converging to zero sequence $\{\varepsilon_k>0\}$ and a sequence $\{K_k\}$ of finite subsets of $G$. There exists a congruent sequence of tilings $\{\tilde{\mathcal{T}_k}\}$ of $G$ such that
	shapes of $\tilde{\mathcal{T}_k}$ are $(K_k,\varepsilon_k)-$invariant.
\end{lemma}

\subsection{Topological entropy for non-compact subsets}

By resembling the definition of Hausdorff dimension, Bowen \cite{bowen1973topological} introduced a definition of topological entropy on subsets for $\mathbb{Z}$-actions. This definition is also known as dimensional entropy and has plenty applications to thermodynamical formulism, fractal geometry, multi-fractal analysis etc. See \cite{barreira2007nonuniform, przytycki2010conformal} for example.

For amenable group actions, Bowen's topological entropy was introduced in \cite{zheng2016bowen} recently.

For  $F\in F(G)$, define
\[
\rho_{F}(x,y) = \max\{\rho(sx,sy) \mid s\in F\}.
\]
Let $Y\subset X$ and $\mathcal{F}=\{F_n\}$ be a F{\o}lner sequence. For $\varepsilon>0$ and $N\in \mathbb{N},$ denote $\mathcal{C}_{N}(Y,\varepsilon,\mathcal{F})$  the collection of all finite or countable covers $\mathcal{C}=\{B_{F_{n_i}}(x,\varepsilon)\}$ of $Y$ with $n_i\geq N.$ For $s>0$, denote
\[
M(Y,\varepsilon,N,s,\mathcal{F}) := \inf\limits_{\mathcal{C}\in \mathcal{C}_{N}(Y,\varepsilon,\mathcal{F})} \sum_{B_{F_m}(x,\varepsilon)\in \mathcal{C}}e^{-s|F_m|}.
\]
The value $M(Y,\varepsilon,N,s,\mathcal{F})$ does not decrease as $N$ increases, hence we can define the following 
\[
M(Y,\varepsilon,s,\mathcal{F}) = \lim\limits_{N\rightarrow \infty}M(Y,\varepsilon,N,s,\mathcal{F}).
\]

It is easy to check there exists a critical value of $s$ such that $M(Y,\varepsilon,s,\mathcal{F})$ jumps from $+\infty$ to $0.$
Let 

\begin{align*}
h_{top}^{B}(Y,\varepsilon,\mathcal{F}) &:= \inf \{s\mid M(Y,\varepsilon,s,\mathcal{F}) = 0 \}\\
&:= \sup \{s\mid M(Y,\varepsilon,s,\mathcal{F}) = \infty \}.
\end{align*}
Clearly $h_{top}^{B}(Y,\varepsilon,\mathcal{F})$ does not decrease as $\varepsilon$ decreases, hence the following limit exists

 \[h_{top}^{B}(Y,\mathcal{F}) = \lim\limits_{\varepsilon \rightarrow 0}h_{top}^{B}(Y,\varepsilon,\mathcal{F}),
\]
and we call it (Bowen) topological entropy of $Y$ with respect to $\mathcal{F}.$

\subsection{Specification}
In this subsection, we will recall the specification property for general group actions, which is from \cite[Section 6]{chung2015homoclinic}.

Let $\alpha$ be a continuous $G$-action on a compact metric space $X$ with metric $\rho.$ The action has the specification property if  for every $\varepsilon>0,$ there is a nonempty finite subset $F=F(\varepsilon)$ of $G$ with the following property: for any finite collection of finite subsets $F_1, F_2,\cdots, F_m$ of $G$ with
\begin{align}
FF_i\cap F_j=\emptyset \qquad 1\leq i\neq j\leq m, \label{specification1}
\end{align}
and for any collection of points $x^1,x^2,\dots,x^m \in X,$ there is a point $y\in X$ with
\begin{align}
\rho(sx^{i},sy)\leq \varepsilon \quad \text{for all}\quad s\in F_i,\ 1\leq i\leq m. \label{con1}
\end{align}

\subsection{Metric Entropy}

Let $\mu$ be an invariant measure and $\beta$ be a finite measurable partition of $X.$ Denote $H_{\mu}(\beta) = -\sum_{B\in \beta}\mu(B)\log \mu(B)$. For $F\in F(G)$, denote $\beta_{F} = \bigvee_{s\in F}s^{-1}\beta$. The metric entropy of $\mu$ with respect to $\beta$ is defined by 
\[
h_{\mu}(\beta, G) = \lim\limits_{n\rightarrow \infty}\frac{1}{|F_n|}H_{\mu}(\beta_{F_n}),
\]
where $\mathcal{F} = \{F_n\}$ is a F{\o}lner sequence and the value $h_{\mu}(\beta,G)$ is independent of the choice of $\mathcal{F}$. The metric entropy of the system $(X,G,\mu)$ is defined as
\[
h_{\mu}(X,G) = \sup_{\beta}h_{\mu}(\beta,G).
\]
The {\itshape entropy map} is a map associating an invariant measure $\mu\in M(X,G)$ to its metric entropy $h_{\mu}(X,G)\in [0,\infty)\cup \{\infty\}$.
For more information about ergodic theory for group actions, readers may see \cite{ollagnier2007ergodic,kerr2016ergodic}.

\subsection{Separated sets}
Let $F\in F(G)$ and $\delta>0, \varepsilon>0$. A subset $\Gamma \subset X$ is $(\delta,F,\varepsilon)-separated$ if for $x\neq y \in \Gamma,$
\[
\frac{|\{s\in F\mid \rho(sx,sy)>\varepsilon\}|}{|F|} \geq \delta.
\]

For $\mu\in M(X,G)$, by $\mathcal{N}(\mu)$ we denote the family of all weak$^*$ neighborhoods of $\mu$.
  Given $C\in\mathcal{N}(\mu)$, we define
\begin{gather}
	X_{F,C} := \{x\in X\mid \mathcal{E}_{F}(x)\in C\},\\
	N(C;F,\varepsilon) := \text{ maximal cardinality of an } (F,\varepsilon)-separated \text{ subset of } X_{F,C}, \\
	N(C;\delta, F,\varepsilon) := \text{ maximal cardinality of a } (\delta, F,\varepsilon)-separated \text{ subset of } X_{F,C}.
\end{gather}

\begin{definition}
	
An $f$-neighborhood of $\mu\in M(X)$ is the set of the form 
\[
F^{(\alpha)} := \{\nu\in M(X) \mid \big|\langle f_i,\mu\rangle - \langle f_i,\nu\rangle\big| \leq \alpha\varepsilon_i\},
\]
where $\alpha>0, \varepsilon_i>0, f_i\in C(X,\mathbb{R}), \ i=1,\dots,k$ and $\norm{f_i}\leq 1$ for each $i$, where $\norm{f_i}=\sup_{x\in X}|f_{i}(x)|$.
\end{definition}
The $f$-neighborhoods form a neighborhood base for the weak$^*$ topology on $M(X)$, which is the topology we use.

The following lemma will be needed in Section 5.
\begin{lemma}\cite[Lemma 2.6]{zhang2018topological}\label{eslemma1}
	
	Let $(X,G)$ be a dynamical system. Let $\mu\in M(X,G), \delta^{*}>0, \varepsilon^{*} >0,\xi>0.$ Let $0< \delta < \min\{\frac{1}{2},\frac{\xi}{3},\frac{\delta^*}{2}\}$, $F\in F(G)$ and $\Gamma \subset X_{F,\mathcal{B}(\mu,\xi)}$ be a $(\delta^*,F,\varepsilon^*)-$separated set. Then for any $F^{\prime} \subset F$ with $\frac{|F^{\prime}|}{|F|}> 1- \delta$, $\Gamma$ is a $(\frac{\delta^*}{2},F^{\prime},\varepsilon^*)-$separated set and $\Gamma \subset X_{F^{\prime},\mathcal{B}(\mu,2\xi)}$.
\end{lemma}
We remark that the statement of Lemma \ref{eslemma1} is a little bit different from \cite[Lemma 2.6]{zhang2018topological}. But the proof of \cite[Lemma 2.6]{zhang2018topological} gives us the above conclusion.

\subsection{Approximation by Ergodic Measures}
\begin{definition}
	The measure $\nu\in M(X,G)$ is entropy-approachable by ergodic measures if for any neighborhood $C\in\mathcal{N}(\nu)$ and each $h^*< h_{\nu}(X,G),$ there exists a measure
	$\mu\in E(X,G)\cap C$ such that $h_{\mu}(X,G) > h^*.$ The ergodic measures are  entropy-dense if each $\nu\in M(X,G)$ is entropy-approachable by ergodic measures.
\end{definition}

For amenable group actions, we prove the following result.
\begin{theorem}\label{approthm}
Let $(X,G)$ be a dynamical system.	Suppose the dynamical system has the specification property. Then the ergodic measures are entropy dense.
\end{theorem}
\begin{proof}
We will prove this theorem in the Appendix.	
\end{proof}

\section{Uniform separation property}

Uniform separation property for $\mathbb{Z}-$actions was introduced by Pfister and Sullivan in \cite{pfister2007topological}.
In this section, we will define the uniform separation property for amenable group actions.
\begin{definition}
	The dynamical system $(X,G)$ has {\itshape uniform separation property} if the following holds. Let $\{K_n\}$ be a tempered F{\o}lner sequence. For any $\eta >0,$ there exists $\varepsilon^*>0$ and $\delta^*>0$ such that for $\mu\in E(X,G)$ and any neighborhood $C\in\mathcal{N}(\mu)$, there exist $n_{C;\mu,\eta}^{*}\in\mathbb{N},$ such that for $n\geq n_{C;\mu,\eta}^{*}$,
	\[
	N(C;\delta^*,K_n,\varepsilon^*) \geq e^{|K_n|(h_{\mu}(X,G) - \eta)}.
	\]
\end{definition}

{\itshape Remark:} Note that uniform separation property implies $h_{top}(X,G) < \infty.$ Indeed for $\mu\in E(X,G)$ we have
\[
 e^{|K_n|(h_{\mu}(X,G) - \eta)} \leq N(C;\delta^*,K_n,\varepsilon^*) \leq N(X;\delta^*,K_n,\varepsilon^*) \leq N(X;K_n,\varepsilon^*),
\]
where $N(X;\delta^*, K_n,\varepsilon^*)$ is the maximal cardinality of a $(\delta^*, K_n,\varepsilon^*)-$separated set of $X$ and $N(X;K_n,\varepsilon^*)$ is the maximal cardinality of a $(K_n,\varepsilon^*)-$separated set of $X$. By \cite[Lemma 9.33]{kerr2016ergodic}, there is an $M>0$ such that for all nonempty finite set $F\subset G$ one has 
\[
N(X;F,\varepsilon^*) \leq M^{|F|}N(X;{\{e\}},\varepsilon^*/2) .
\]
Thus
\[
h_{\mu}(X,G) \leq \limsup_{n\rightarrow \infty}\frac{\log  N(X;K_n,\varepsilon^*)}{|K_n|} + \eta  \leq M + \frac{\log N(X;\{e\},\varepsilon^*/2)}{|K_n|} + \eta< \infty.
\]

The following lemma appears as \cite[Lemma 1.5.4]{shields1996ergodic}.
\begin{lemma}\label{combination}
	If $\binom{n}{k}$ denotes the number of combinations of $n$ objects taken $k$ at a time and $\delta < 1/2$ then
	\begin{displaymath} 
	\sum_{k\leq \delta n}\binom{n}{k} \leq e^{n\phi(\delta)},
	\end{displaymath}
where $\phi(\delta) = -\delta\log \delta - (1-\delta)\log (1-\delta)$.
\end{lemma}

\begin{theorem}\label{expansive}
	Suppose the action $(X,G)$ is expansive. Then the action has uniform separation property.

\end{theorem}

\begin{proof}
Let $\{K_n\}$ be a tempered F{\o}lner sequence. Let $\tau^*$ be the expansive constant. Then for any finite Borel partition $\mathcal{A}$ with $diam(\mathcal{A})<\tau^*$, $h_{\mu}(X,G) = h_{\mu}(\mathcal{A},G)\ \forall \mu\in M(X,G).$ Let $\eta > \eta^\prime>0.$ The value of $\eta^\prime$ will be fixed in (\ref{choiceofeta}). Let $\nu\in M(X,G)$. We first construct a neighborhood $W_{\nu}\subset M(X)$  of $\nu.$ Since $M(X,G)$ is compact, we can cover $M(X,G)$ by finite neighborhoods $W_{\nu_1},W_{\nu_2},\dots, W_{\nu_n}$. It is enough to prove the result for an ergodic measure $\mu\in W_{\nu}.$

Let $\mathcal{A}=\{A_1,A_2,\dots,A_k\}$ be a finite Borel partition with $diam(\mathcal{A}) < \frac{\tau^*}{2}.$

Choose $\delta^*$ small such that $2\eta^\prime + \delta^* < \frac{1}{2}$ and 
\begin{align}
\phi(2\eta^\prime + \delta^*) + (2\eta^\prime + \delta^*)\log(2k-1) < \eta - \eta^\prime, \label{choiceofeta}
\end{align}
where $\phi$ is as described  in Lemma \ref{combination}.

Since $\nu$ is regular, for every $j=1,\cdots,k$ we can find a compact set $V_j\subset A_j$ with $\nu(A_j\setminus V_j) < \frac{\eta^\prime}{4k\log (2k)}.$ Define $\varepsilon^* = \min \{dist(V_i,V_j)\mid 1\leq i< j\leq k\}/2.$ There exists $n^*$, such that for $n\geq n^*$ we have
\begin{align}
\frac{\eta^{\prime}}{4\log k} \geq \frac{\log 2}{|K_n|\log 2k}. \label{usest1} 
\end{align}

For each $j=1,\cdots,k$ we pick an open neighborhood $U_j$ of $V_j$ with $diam(U_j)<\tau^*$. We can do this in such a way that if $x\in U_i$ and $y\in V_j$ for some $i\neq j$ then $\rho(x,y)>\varepsilon^*$. Let $K = X\setminus \cup_{j=1}^{k}U_j$ and $\mathcal{A}^{\prime}$ a Borel partition including all $U_i, i=1,\dots,k$ and the non-empty intersection $A_i\cap K$. Then $K$ is a closed subset of $X$ and $\nu(K) < \frac{\eta^\prime}{4\log (2k)}.$ The indicator function $I_K$ is upper semi-continuous. We define the neighborhood $W_{\nu}$ of $\nu$ by 
\[
W_{\nu} := \left\{m\in M(X) : \int I_K\dd m \leq \int I_K\dd\nu + \frac{\eta^{\prime}}{4\log(2k)}\right\}.
\]

By the construction, $\mathcal{A}^{\prime}$ is a partition with no more that $2k$ elements and $diam(\mathcal{A}^\prime) < \tau^*.$ For convenience, we will label each element in $\mathcal{A}^{\prime}_{K_n}$ by a word $w$ of length $K_n$ over an alphabet $A$ of at most $2k$ letters. The letters $1,2,\dots,k$ label $U_1, U_2,\dots, U_k$ and the other letters label the non-empty atoms among $A_1\cap K,\dots,A_k\cap K.$ We will define a map $\varPhi: X \rightarrow A^{K_n}$ by
\[
\varPhi(x)_{s} := j \quad \text{ if } sx \text{ is in the atom labeled by } j.
\]
Since $diam(\mathcal{A}^{\prime}) < \tau^*,$ we will choose $n^*$ such that for $n\geq n^*$
\[
H_{\mu}(\mathcal{A}^{\prime}_{K_n}) \geq |K_n|(h_{\mu}(X,G) - \eta^{\prime}/2).
\]
By definition of $W_{\nu}$ we have
\[
\mu(K) \leq \nu(K) + \frac{\eta^{\prime}}{4\log(2k)} \leq \frac{\eta^{\prime}}{2\log (2k)}.
\]

Let $Y_n = \left\{x\in X : |\{s\in K_n : \varPhi(x)_s > k\}|\leq \eta^{\prime}|K_n|\right\}.$ By the pointwise  ergodic theorem, $\mu(Y_n) \rightarrow 1.$ Let $C\in\mathcal{N}(\mu)$. There is $n^*$ such that for $n\geq n^*$ we have 
\begin{equation}\mu(X_{K_n,C}) > 1 - \frac{\eta^{\prime}}{3\log(2k)}\end{equation}
 and
  \begin{equation}
\mu(X\setminus(X_{K_n,C}\cap Y_n)) < \frac{\eta^{\prime}}{2\log(2k)} - \frac{\log 2}{|K_n|\log(2k)}.  \label{equationneighborhood}
\end{equation}

Set $A_{0}^{n} =X\setminus(X_{K_n,C}\cap Y_n).$ By conditioning with respect to the partition $\{A_{0}^{n},X\setminus A_{0}^{n}\}$, we obtain
\begin{align*}
H_{\mu}(\mathcal{A}^{\prime}_{K_n}) &\leq \log 2 + \mu(A_{0}^{n})H_{\mu(\cdot|A_{0}^{n})}(\mathcal{A}^{\prime}_{K_n}) + \mu(X\setminus A_{0}^{n})H_{\mu(\cdot|X\setminus A_{0}^{n})}(\mathcal{A}^{\prime}_{K_n}).
\end{align*}

Since the number of atoms in $\mathcal{A}^\prime$ is at most $2k$ we have $\mu(A_{0}^{n})H_{\mu(\cdot|A_{0}^{n})}(\mathcal{A}^{\prime}_{K_n}) \leq  \big(\frac{\eta^{\prime}}{2\log(2k)} - \frac{\log 2}{|K_n|\log(2k)}\big)|K_n|\log(2k)\leq \frac{\eta^{\prime}}{2}|K_n| - \log 2.$ Thus for $n\geq n^*,$ 
\[
\mu(X\setminus A_{0}^{n})H_{\mu(\cdot|X\setminus A_{0}^{n})}(\mathcal{A}^{\prime}_{K_n}) \geq |K_n|(h_{\mu}(X,G) - \eta^{\prime}).
\]

Let $\varPhi_n$ denote the image of $X_{K_n,C}\cap Y_n$ by the map $\varPhi.$ Since $\log |\{A\in \mathcal{A}^{\prime}_{K_n} : A\cap (X\setminus A_{0}^n) \neq \emptyset\}| \geq H_{\mu(\cdot|X\setminus A_{0}^{n})}(\mathcal{A}^{\prime}_{K_n})$, we have 

\begin{align}
|\varPhi_n| \geq e^{|K_n|(h_{\mu}(X,G) - \eta^{\prime})}. \label{uses3}
\end{align}

A word on $K_n$ over a finite alphabet $\mathcal{A}$ is a function $w: K_n \rightarrow \mathcal{A}$. Given two different words $w, w^{\prime}$ on a common domain  and a common alphabet, we say $d_{K_n}^{H}(w, w^{\prime})$ to be the number of entries $q$ in $K_n$ such that $w(q)\neq w^{\prime}(q)$.
 Let $\varPhi^{\prime}_n \subset  \varPhi_n$ of maximal cardinality such that $d_{K_n}^{H}(w,w^{\prime}) \geq (2\eta^{\prime}+\delta^*)|K_n|$ for any $w\neq w^{\prime} \in \varPhi^{\prime}_n.$ By Lemma \ref{combination} and the choice of $\eta^{\prime}\text{ and } \delta^*$ we have
 \begin{align*}
 	|\varPhi^{\prime}_n| &\geq \frac{|\varPhi_n|}{\sum_{j\leq (2\eta^{\prime} + \delta^*) |K_n|}\binom{|K_n|}{j}(2k-1)^{(2\eta^{\prime} + \delta^*)|K_n|}}\\
 	 &\geq e^{|K_n|\big(h_{\mu}(X,G) - \eta^{\prime} - \phi(2\eta^{\prime} + \delta^*)- (2\eta^{\prime} + \delta^*)\log(2k-1)\big)} \geq e^{|K_n|(h_{\mu}(X,G)-\eta)}. 
\end{align*} 
 Let $\Gamma_n$ be defined by selecting exactly one point of $X_{K_n, C} \cap Y_n$ from each atom of $\mathcal{A}_{K_n}^{\prime}$ labeled by a word in $\varPhi_n^{\prime}.$ Then $\Gamma_n$ is $(\delta^*,K_n,\varepsilon^*)-$separated and
\[
|\Gamma_n| \geq e^{|K_n|(h_{\mu}(X,G) - \eta)}.
\]
The proof is finished.
\end{proof}

\begin{proposition}\label{separatinges1}
	Let $(X,G)$ be a dynamical system. Assume the system has uniform separation property and that the ergodic measures are entropy dense. Let $\{K_n\}$ be a tempered F{\o}lner sequence. For any $\eta >0$, there exist $\delta^*>0$ and $\varepsilon^*>0$ so that for $\mu\in M(X,G)$ and any neighborhood $C\in \mathcal{N}(\mu)$, there exists $n^{*}_{C;\mu,\eta}$ such that for $n\geq n^{*}_{C;\mu,\eta}$
	\[
	N(C;\delta^*,K_n,\varepsilon^*) \geq e^{|K_n|(h_{\mu}(X,G) - \eta)}.
	\]
	Furthermore, for any $\mu\in M(X,G)$ we have
	\[
	h_{\mu}(X,G) \leq \lim\limits_{\varepsilon\rightarrow 0}\lim\limits_{\delta\rightarrow 0}\inf_{C\in\mathcal{N}(\mu)}\liminf\limits_{n\rightarrow \infty}\frac{1}{|K_n|}\log N(C;\delta,K_n,\varepsilon).
	\]
\end{proposition}

\begin{proof}
	Let $\eta>0$ and $\mu \in C$. If $\mu$ is ergodic, then the statement is true by the definition of  uniform separation property. If $\mu$ is not ergodic, then choose an ergodic $\nu \in C$ and $h_{\mu}(X,G) < h_{\nu}(X,G) + \frac{\eta}{2}.$  We can just choose $n_{C;\mu,\eta}^* = n_{C;\nu,\frac{\eta}{2}}^*.$ The second statement is a consequence of the first statement.
\end{proof}

Next theorem shows one example which has the specification and uniform separation properties.
\begin{theorem}
	Let $\Gamma$ be a  countable discrete group and $f$ an element
	of $\mathbb{Z}\Gamma$ invertible in $l^{1}(\Gamma, \mathbb{R}).$  Then the action of $\Gamma$ on $X_f$ which is the Pontryagin dual of $\mathbb{Z}\Gamma/\mathbb{Z}\Gamma f$ has
	the specification and uniform separation properties.
\end{theorem}

\begin{proof}
	By \cite[Theorem 1.2]{ren2018}, the system has the specification property and also it is expansive. Then by Theorem \ref{expansive}, the system has uniform separation property.
\end{proof}

\subsection{The entropy map}
In this part we will study the entropy map when uniform separation property holds. Let $\mathcal{F}=\{F_n\}$ be a F{\o}lner sequence. For $\varepsilon >0, \delta>0$ and $\nu\in M(X,G)$, denote
\begin{equation}
	\underline{s}(\nu;\delta,\varepsilon,\mathcal{F}) = \inf_{C\in\mathcal{N}(\nu)}\liminf_{n\rightarrow \infty}\frac{1}{|F_n|}\log N(C;\delta,F_n,\varepsilon) \label{dfnse1}	
\end{equation}
and 
\begin{equation}
	\underline{s}(\nu;\mathcal{F}) := \lim\limits_{\varepsilon \rightarrow 0}\lim\limits_{\delta \rightarrow 0} \underline{s}(\nu;\delta,\varepsilon,\{F_n\}). \label{dfnse2}	
\end{equation}

We define $\overline{s}(\nu;\delta,\varepsilon,\mathcal{F})$ and $\overline{s}(\nu;\mathcal{F})$ by taking $\limsup$ instead of $\liminf$ of \eqref{dfnse1} and \eqref{dfnse2}. If $\underline{s}(\nu;\mathcal{F})= \overline{s}(\nu;\mathcal{F})$, then denote $s(\nu;\mathcal{F})$  the common value.

	Let $(X,G)$ be a topological dynamical system and $\mu\in M(X,G)$. Let $\mathcal{F}=\{F_n\}$ be a F{\o}lner sequence. Let $E_n$ be a sequence of $(F_n, \varepsilon)-$separated subsets and define
	$$\nu_{n} := \frac{1}{|F_n||E_n|}\sum_{x\in E_n}\sum_{s\in F_n}\delta_{sx}.$$
Assume $\nu_n \rightarrow \mu.$ Then
	$$\limsup_{n\rightarrow \infty.}\frac{1}{|F_n|}\log|E_n| \leq h_{\mu}(X,G).$$

For $\mu\in M(X,G)$ and a F{\o}lner sequence $\mathcal{F}$, denote
\[
\overline{s}^{\prime}(\nu;\varepsilon,\mathcal{F}) = \inf_{C\in \mathcal{N}(\nu)}\limsup_{n\rightarrow \infty}\frac{1}{|F_n|}\log N(C;F_n,\varepsilon)
\]
 and
 \[
 \overline{s}^{\prime}(\nu;\mathcal{F}) =\lim\limits_{\varepsilon \rightarrow 0}\overline{s}^{\prime}(\nu;\varepsilon,\mathcal{F}) .
 \]
\begin{proposition}\label{sec3.1-prop}
	Let $(X,G)$ be a topological dynamical system and $\mu \in M(X,G)$. Then for any F{\o}lner sequence $\mathcal{F}=\{F_n\}$ we have
	$$\overline{s}^{\prime}(\mu; \mathcal{F})\leq h_{\mu}(X,G).$$
\end{proposition}

\begin{proof}
	The proof is similar to the proof of \cite[Proposition 3.1]{pfister2007topological}. 
	
	If $h_{\mu}(X,G) = \infty,$ then there is nothing to prove. Let $h_{\mu}(X,G) < \infty.$ Suppose that
	$$\lim\limits_{\varepsilon \rightarrow 0}\inf_{C \in \mathcal{N}(\mu)}\limsup_{n\rightarrow \infty}\frac{1}{|F_n|}\log N(C;F_n,\varepsilon) > h_{\mu}(X,G).$$
There exist $\varepsilon^{*}>0$ and $\eta >0$ such that for $0<\varepsilon \leq \varepsilon^{*},$
	$$\inf_{C \in \mathcal{N}(\mu)}\limsup_{n\rightarrow \infty}\frac{1}{|F_n|}\log N(C; F_n, \varepsilon) > h_{\mu}(X,G) + 2 \eta.$$
	Let $0 < \varepsilon < \varepsilon^{*}$. There exists a decreasing sequence of convex closed neighborhoods $\{C_n\}$ of $\mu$ such that
	$\bigcap_{n}C_n = \{\mu\}$ and 
	\begin{equation}\label{sec3.1-prop1}
	\limsup_{n\rightarrow \infty}\frac{1}{|F_n|}\log N(C_n; F_n, \varepsilon) \geq h_{\mu}(X,G) + 2\eta.	
	\end{equation}
Let $E_n$ be a $(F_n, \varepsilon)-$separated set of $X_{F_n,C}$ with maximal cardinality, and define
	
	$$\nu_{n} := \frac{1}{|F_n||E_n|}\sum_{x\in E_n}\sum_{s\in F_n} \delta_{sx} \in C_n.$$
By the choice of $\{C_n\}$, we have $\lim\limits_{n\rightarrow \infty} \nu_{n} = \mu.$ Using the standard arguments in the proof of the variational principle(for example see \cite[Page 227]{kerr2016ergodic}), we have
	$$\limsup_{n\rightarrow \infty}\frac{\log |E_n|}{|F_n|} = \limsup_{n\rightarrow \infty}\frac{1}{|F_n|}\log N(C_n; F_n, \varepsilon) \leq h_{\mu}(X,G),$$
	which contradicts (\ref{sec3.1-prop1}).
\end{proof}

\begin{proposition}\label{sepes1}
	
	Let $\mathcal{F}=\{F_n\}$ be a  F{\o}lner sequence and $\mu\in E(X,G)$. Then for every $h^*<h_{\mu}(X,G)$, there exist $\delta^*>0, \varepsilon^*>0$ such that for any neighborhood $C\in\mathcal{N}(\mu)$, there exists $n_{C}^*$ such that for any $n\geq n_{C}^*$ there exists a $(\delta^*,F_n,\varepsilon^*)-$separated set $\Gamma_n\subset X_{F_n,C}$ satisfying 
	$|\Gamma_n|\geq e^{h^*|F_n|}$.
\end{proposition}

\begin{proof}
The proof is an adaptation of the proofs \cite[Theorem 8.6]{walters2000introduction} and \cite[Proposition 2.1]{pfister2004large} to the setting of countable amenable group actions.

By the mean ergodic theorem, for $n$ large, $X_{F_n,C} \neq \emptyset.$ If $h^* \leq 0,$ we may just take $\Gamma_n$ to be a one-point set.

Let $h_{\mu}(X,G) >0$ and $0<h^*<h_{\mu}(X,G).$ Choose $h^{\prime}, h^{\prime\prime},h_1,h_2 $ satisfying $h^* < h^{\prime} < h^{\prime\prime} < h_1 <h_2< h_{\mu}(X,G).$  Take an $f-$neighborhood $F^{(1)} \subset C$ of $\mu$ corresponding to $\{f_i,\varepsilon_i : i=1,2,\dots,l\}.$

Take a partition $\gamma=\{A_1,\dots,A_k\}$ such that $h_{\mu}(\gamma,G) >  h_2$. Choose $\delta^* >0$ satisfying
\begin{align}
\phi(2\delta^*) + 2\delta^{*}{(k+1)} < h^{\prime} - h^{*}.\label{delta}
\end{align}

Since $\mu$ is regular, for any $\theta>0$ and for every $1\leq j\leq k$, there exists a compact set $B_j \subset A_j$ with $\mu(A_j\setminus B_j) <\theta.$ The value of $\theta$ will be fixed later. Let $\beta = \{B_0,B_1,\dots,B_k\}$, where $B_0 = X\setminus \cup_{j=1}^{k}B_j$. By \cite[Proposition 9.4]{kerr2016ergodic},
\[
h_{\mu}(\gamma,G) \leq h_{\mu}(\beta,G) + H(\gamma|\beta),
\]
where $H(\gamma|\beta) = - \sum_{B\in \beta}\sum_{A\in\gamma}\mu(A\cap B)\log\frac{\mu(A\cap B)}{\mu(B)}$ is the conditional entropy. Hence we assume $\theta>0$ to be so small that $H(\gamma|\beta) < h_2 - h_1$. By this choice of $\theta$, we obtain $h_{\mu}(\beta,G)> h_1$. We also assume $\theta < \min\{\frac{\delta^*}{4k},\frac{h^{\prime\prime} - h^{\prime}}{12(k+1)\log (k+1)}\}$.

Define $\varepsilon^{*} = \min\{dist(B_i,B_j) \mid 1\leq i< j\leq k \}/2.$  

 Fix $\psi \in L^{1}(X,\mathbb{R})$. For a finite set $F\subset G$, define $S_{F}\psi(x)=\sum_{s\in F}\psi(sx)$. By the mean ergodic theorem, for every $\sigma >0$ it holds
\begin{align}
\lim\limits_{n\rightarrow \infty} \mu\{x\in X: \big|\frac{1}{|F_n|}S_{F_n}\psi(sx) - \int \psi \dd\mu\big|> \sigma\} = 0. \label{mean}
\end{align}
Then there exists $n_1^*\in \mathbb{N}$ such that for $n>n_1^*$,
\begin{align}
	\mu(X_{F_n,F^{(1)}}) > 1- (h^{\prime\prime} - h^{\prime})/(3\log (k+1)). \label{displa1}
\end{align}

Take $\sigma = \min\{\frac{\delta^*}{4},\frac{h^{\prime\prime} - h^{\prime}}{12(k+1)\log (k+1)}\}$. 
Define $\varphi(x) := I_{\cup_{j=1}^{k}B_j}(x)$
where we write $I_{A}$ for the indicator function on $A$.  By (\ref{mean}), there exists $n_2^{*}\in\mathbb{N}$ such that for $n\geq n_2^{*}$ there exists a measurable subset $X_n^{\prime} $ with $\mu(X_n^{\prime} ) > 1- (h^{\prime\prime} - h^{\prime})/(3\log (k+1))$ and for $x\in X_n^{\prime}$
we have
\[
\big|\frac{1}{|F_n|}S_{F_n}\varphi(x) - \mu(\cup_{j=1}^{k}B_j)\big| \leq \sigma.
\]
Let $n^*$ be so large that for $n \geq n^*$ it holds
\[
H_{\mu}(\beta_{F_n}) > |F_n|h^{\prime\prime} \text{\  and \ } \frac{h^{\prime\prime} - h^{\prime}}{6(k+1)\log (k+1)} > \frac{\log 2}{|F_n|\log (k+1)} \text{ and }\delta^{*}|F_n| > 2.
\]
Pick $n_{C}^* > \max\{n_1^*,n_2^*,n^*\}$ such that for $n\geq n_{C}^*$ we have
\[
\mu(X\setminus X_{F_n,F^{(1)}}) \leq (h^{\prime\prime} - h^{\prime})/(3\log (k+1)),\ \mu(X\setminus X_n^{\prime}) \leq (h^{\prime\prime} - h^{\prime})/(3\log (k+1))
\]
and for any $x\in X_n^{\prime}$
\begin{align}
\frac{1}{|F_n|}S_{F_n}\varphi(x)>  1- k\theta-\sigma. \label{xes}
\end{align}

Let $A_{0}^{n} = X\setminus(X_{F_n,F^{(1)}}\cap X_n^{\prime})$. Then for $n\geq n_{C}^*$,
\[
 \mu(A_{0}^n)  <\frac{2(h^{\prime\prime} - h^{\prime})}{3\log (k+1)} <\frac{(h^{\prime\prime} - h^{\prime})}{\log (k+1)} - \frac{\log 2}{|F_n|\log (k+1)}.
\]

By considering  the partition $\{A_{0}^{n},X\setminus A_{0}^{n}\}$, we obtain
\begin{align*}
H_{\mu}(\beta_{F_n}) &\leq \log 2 + \mu(A_{0}^{n}) H_{\mu(\cdot|A_{0}^{n})}(\beta_{F_n}) + \mu(X\setminus A_{0}^{n})H_{\mu(\cdot|X\setminus A_{0}^{n})}(\beta_{F_n})\\
&\leq \log 2 + \Big(\frac{h^{\prime\prime} - h^{\prime}}{\log (k+1)} - \frac{\log 2}{|F_n|\log (k+1)}\Big)|F_n|\log (k+1) + H_{\mu(\cdot|X\setminus A_{0}^{n})}(\beta_{F_n})\\
&\leq (h^{\prime\prime} - h^{\prime})|F_n| + H_{\mu(\cdot|X\setminus A_{0}^{n})}(\beta_{F_n}).
\end{align*}
Then we have $H_{\mu(\cdot|X\setminus A_{0}^{n})}(\beta_{F_n}) \geq h^{\prime}|F_n|$, which implies that the number of elements of the set
\begin{align}
\beta_n :=\{B\in \beta_{F_n}\mid B\setminus A_{0}^{n} \neq \emptyset\} \label{partition} 
\end{align}
is at least $e^{h^{\prime}|F_n|}$.

Define a map $\varPhi : X \rightarrow \{0,1,\dots,k\}^{F_n}$ by
\[
\varPhi(x)_{s} := j \text{ if } sx\in B_j,\text{ for } j=0,1,\dots,k. 
\]
Denote by $\varPhi_n$  the image of $X_{F_n,F^{(1)}}\cap X_n^{\prime} .$ Then by (\ref{partition}) it holds
\[
|\varPhi_n| \geq e^{h^{\prime}|F_n|}.
\]
Recall that the  Hamming distance between $w=\varPhi(x)$ and $w^{\prime}=\varPhi(x^{\prime})$ is the number of entries $s\in F_n$ where $w$ and $w^{\prime}$ differ. Let 
$\varPhi_n^{\prime}\subset \varPhi_n$ be the subset of maximal cardinality with
$d_{F_n}^{H}(w,w^{\prime}) \geq 2\delta^{*}|F_n|.$
By Lemma \ref{combination} and the choice of $\delta^*$ we have
 \[
 |\varPhi^{\prime}_n| \geq \frac{e^{|F_n|h^{\prime}}}{e^{\phi(2\delta^*)|F_n|}k^{2\delta^{*}|F_n|}} \geq e^{h^{*}|F_n|}.
 \]
 
 Let $\Gamma_n$ be defined by selecting exactly one point of $X_{F_n,F^{(1)}}\cap X_n^{\prime} $ from each atom of $\beta_{F_n}$ labeled by a word in $\varPhi_{n}^{\prime}.$  For any $x\in \Gamma_n$ it holds  $\big|\{s\in F_n \mid \varPhi(x)_s = 0\} \big| \leq (k\theta+\sigma)|F_n|$ and for different $x,y\in \Gamma_n$ we have $|\{s\in F_n \mid \varPhi(x)_s \neq \varPhi(y)_s \}| \geq 2\delta^*|F_n|$. Hence  for $x\neq y\in \Gamma_n$ we obtain
 \[
 |\{s\in F_n \mid \varPhi(x)_s \neq \varPhi(y)_s \in \{1,2,\dots,k \}| \geq (2\delta^* - 2(k\theta+\sigma))|F_n| \geq \delta^*|F_n|,
 \]
 which means that $x,y$ are $(\delta^*,F_n,\varepsilon^*)-$separated.
 Hence $\Gamma_n$ is a $(\delta^*,F_n,\varepsilon^*)-$separated set with $|\Gamma_n| \geq e^{h^{*}|F_n|}.$ 
\end{proof}

\begin{corollary} \label{entropy1}
	Let $(X,G)$ be a topological dynamical system and $\mathcal{F}$ be a F{\o}lner sequence. For $\nu \in E(X,G),$ 
	\begin{align*}
	h_{\nu}(X,G) = h_{\nu}(X,\mathcal{F}) &=\underline{s}(\nu;\mathcal{F})\\
	& = \overline{s}(\nu;\{F_n\}).
	\end{align*}
\end{corollary}

\begin{proof}
	 Using the fact the inequality $N(C;\delta,F,\varepsilon) \leq N(C;F,\varepsilon)$ holds for every $F\in F(G)$, we have
	\[
	\overline{s}(\nu;\delta,\varepsilon,\mathcal{F}) \leq \overline{s}^{\prime}(\nu;\varepsilon,\mathcal{F}).
	\]
Using Proposition \ref{sec3.1-prop} we conclude that
\[
\overline{s}(\nu;\mathcal{F}) \leq \overline{s}^{\prime}(\nu;\mathcal{F}) \leq h_{\mu}(X,G).
\]
By Proposition \ref{sepes1} we obtain
	\begin{align*}
	h_\nu(X,G) \leq \lim\limits_{\varepsilon\rightarrow 0}\lim\limits_{\delta\rightarrow 0}\inf_{C\in \mathcal{N}(\nu)}\liminf_{n\rightarrow \infty}\frac{1}{|F_n|}\log N(C;\delta,F_n,\varepsilon). \label{tempered1}
		\end{align*}
The proof is finished.
\end{proof}
\begin{proposition}
	Let $(X,G)$ be a topological dynamical system. If  uniform separation property condition is true and the ergodic measures are entropy dense, then for a tempered F{\o}lner sequence $\mathcal{K}=\{K_n\}$, $s(\mu,\mathcal{K})$ is well-defined, and $s(\mu,\mathcal{K}) = h_{\mu}(X,\mathcal{K}) = h_{\mu}(X,G)$, for all $\mu\in M(X,G).$
\end{proposition}
\begin{proof}
	Same proof as the proof of Corollary \ref{entropy1}, using Proposition \ref{separatinges1} instead of Proposition \ref{sepes1}.
\end{proof}

Theorem \ref{expansive} tells us that uniform separation property is weaker than expansiveness. The entropy map for expansive amenable group actions is upper semi-continuous (for example see \cite{ren2016local}). The following proposition shows the upper-semicontinuity of the entropy map for systems with uniform separation property.
\begin{proposition}
	Let $(X,G)$ be a topological dynamical system. If uniform separation property condition is true and the ergodic measures are entropy dense, then the entropy map
	\[
	\mu \mapsto h_{\mu}(X,G)
	\]
is upper semi-continuous on $M(X,G)$.
\end{proposition}

\begin{proof}
Let $\{K_n\}$ be a tempered F{\o}lner sequence and $C\in\mathcal{N}(\mu)$. Given $\eta>0,$ by Proposition \ref{separatinges1}, there exists $\delta^*>0$ and $\varepsilon^*>0$  such that 
	\[
	\limsup_{n\rightarrow \infty}\frac{\log N(C;\delta^*,K_n,\varepsilon^*)}{|K_n|} + \eta \geq \sup_{\nu\in C}h_{\nu}(X,G).
	\]
Hence
\[
\inf_{C\in\mathcal{N}(\mu)}\sup_{\nu\in C}h_{\nu}(X,G) \leq \inf_{C\in \mathcal{N}(\mu)}\limsup_{n\rightarrow \infty}\frac{1}{|K_n|}\log N(C;\delta^*,K_n,\varepsilon^*) + \eta.
\]
Since $\eta$ is arbitrary, by Proposition \ref{sec3.1-prop} we obtain
\[
\inf_{C\in\mathcal{N}(\mu)}\sup_{\nu\in C}h_{\nu}(X,G) \leq h_{\mu}(X,G).
\]
The proof is finished.
\end{proof}

\section{Upper bound for $h^B_{top}(G_K(\mathcal{F}),\mathcal{F})$}

\begin{proposition}\label{lemmaforthm1}
	Let $\mathcal{F}=\{F_n\}$ be a F{\o}lner sequence with $\frac{|F_n|}{\log n}\rightarrow \infty$ and $K \subset M(X,G)$ be a closed subset. Let  
			\[^{K}G = \{ x\in X\mid \{\mathcal{E}_{F_n}(x)\} \text{ has a limit point in } K\},\]
			then			
			\[
			h_{top}^{B}(^{K}G, \mathcal{F})\leq \sup\{h_{\mu}(X,G) \mid \mu\in K\}.
			\]
\end{proposition}

\begin{proof}
	Let $ s:= \sup \{h_{\mu}(X,G) \mid \mu\in K\}.$ If $ s= \infty,$ there is nothing to prove. Assume $s < \infty;$ let $s^{\prime} - s = 2\eta.$ Since 
	$N(C; F_n, \varepsilon)$ is a non-increasing function of $\varepsilon,$ by Proposition \ref{sec3.1-prop}, for every $\mu\in M(X,G)$ and $\varepsilon>0$ we have
	$$\inf_{C\in\mathcal{N}(\mu)}\limsup_{n\rightarrow \infty}\frac{1}{|F_n|}\log N(C; F_n, \varepsilon) \leq h_{\mu}(X,G) \text{ for all } \varepsilon>0.$$
	Hence for any $\varepsilon >0$, there exists a neighborhood $C(\mu, \varepsilon)$ of $\mu$ and $n(C(\mu, \varepsilon))\in \mathbb{N}$, such that
	\begin{equation}
	\frac{1}{|F_m|}\log N(C(\mu, \varepsilon);F_m,\varepsilon) \leq h_{\mu}(X,G) + \eta \text{ for all } m \geq n(C(\mu, \varepsilon)).
	\end{equation}
	
	Since a maximal $(F_m,\varepsilon)-$separated set of some $A\subset X$ is also a $(F_m, \varepsilon)-$spanning set of $A$, for any $m \geq n(C(\mu, \varepsilon))$,
	\begin{align}\label{estimatesmall1}
	M(X_{F_m,C(\mu, \varepsilon)},\varepsilon, m, s^{\prime},\{F_n\}) \leq N(C(\mu, \varepsilon); F_m, \varepsilon)e^{-s^{\prime}|F_m|} \leq e^{-\eta |F_m|}.
	\end{align}
	
Since $K$ is compact, given a fixed $\varepsilon,$ we can find a finite open cover of $K$ by sets of the form $C(\mu, \varepsilon),$ say $\mu_1, \mu_2, \dots, \mu_{m_\varepsilon,}$ with $\mu_i \in K.$ If $\{\mathcal{E}_{F_n}(x)\}$ has a limit point in $K$, then for any  $N\in \mathbb{N}$, $x$ is an element of
\[
A_N = \bigcup_{m \geq N}\bigcup_{j=1}^{m_\varepsilon} X_{F_m, C(\mu_j, \varepsilon)}.
\]
Thus for $N \geq \max_{j}N(C(\mu_j,\varepsilon))$, by (\ref{estimatesmall1}), we have,
\begin{align*}
M(^{K}G, \varepsilon, N,s^{\prime}, \mathcal{F}) &\leq \sum_{m\geq N}M(\bigcup_{j=1}^{m_\varepsilon} X_{F_m, C(\mu_j, \varepsilon)},\varepsilon,m,s^{\prime},\mathcal{F})\\
&\leq m_{\varepsilon}\sum_{m\geq N}e^{-\eta|F_m|} < \infty
\end{align*}
which implies that 
\[
h^B_{top}(^{K}G,\varepsilon,\mathcal{F}) \leq s.
\]
The proof is finished.
\end{proof}

\begin{theorem}\label{upper}
	Let $\mathcal{F}=\{F_n\}$ be a F{\o}lner sequence with $\frac{|F_n|}{\log n}\rightarrow \infty$ and $K\subset M(X,G)$ be non-empty compact. Then
		\[ 
		h^B_{top}(G_K(\mathcal{F}),\mathcal{F}) \leq \inf\{h_{\mu}(X,G) : \mu\in K\}.
		 \]
\end{theorem}

\begin{proof}
For all $\mu\in K,$ note that $G_K \subset ^{\{\mu\}}\!\!G.$ By Proposition \ref{lemmaforthm1}
\[
h^B_{top}(G_K(\mathcal{F}),\mathcal{F}) \leq h^{B}_{top}(^{\{\mu\}}G,\mathcal{F})\leq h_{\mu}(X,G)\  \forall \mu\in K,
\]
which means $h^B_{top}(G_K(\mathcal{F}),\mathcal{F}) \leq \inf\{h_{\mu}(X,G) \mid \mu \in K\}.$
\end{proof}

\section{Lower Bound for $h^B_{top}(G_K(\mathcal{F}),\mathcal{F}$)}

\begin{theorem}\label{lowerbound}
	Let $(X,G)$ be a dynamical system with uniform separation property and the specification property. Let $\mathcal{F}=\{F_n\}$ be a F{\o}lner sequence and $K$ be a connected non-empty  closed subset of $M(X,G).$ Then
	\[
	  h^B_{top}(G_K(\mathcal{F}),\mathcal{F}) \geq \inf\{h_{\mu}(X,G) \mid \mu\in K\}. \label{lowerboundes1}
	\]
\end{theorem}
We will show that for any $0< h^* < \inf\{h_{\mu}(X,G) \mid \mu\in K\}$ it holds
\[ 
h^B_{top}(G_K(\mathcal{F}),\mathcal{F}) \geq h^*.
\]
To get this, we will construct a closed subset  $Y_K\subset G_K(\mathcal{F})$ and show that
$ 
h^B_{top}(Y_K,\mathcal{F}) \geq h^*.
$

\subsection{Construction of $Y_K$}

For each $\varepsilon>0$, there exist a finite sequence $\alpha_1,\dots,\alpha_n$ in $K$ such that each point in $K$ is $\varepsilon$ close to some $\alpha_i.$ As $K$ is connected, by repeating some $\alpha_i$, we can choose this so that each point in $K$ is within $\varepsilon$ of some $\alpha_i$ and $D(\alpha_j,\alpha_{j+1}) < \varepsilon$ for each $j.$ Extending this argument, we deduce that there exists a sequence $\{\alpha_j : j =1,2,\dots\}$ in $K$ so that for each $n$ the closure of
$\{\alpha_j : j>n\}$  equals $K$ and 
\[
\lim\limits_{j\rightarrow\infty}D(\alpha_j,\alpha_{j+1}) = 0.
\]
We will construct $Y_K$ such that for each $y\in Y_K$ the sequence \{$\mathcal{E}_{F_n}(y)$\} has the same limit-point set as $\{\alpha_j : j>n\}$ and $h_{top}^B(Y_K,\{F_n\}) \geq h^*.$

Fix a  tempered F{\o}lner sequence $\{K_n\}$. Let $\eta = \inf\{h_{\mu}(X,G) \mid \mu\in K\} - h^{*}$. By Proposition \ref{separatinges1}, we can find $\delta^*>0$ and $\varepsilon^*>0$ so that for every neighborhood $C\in \mathcal{N}(\mu)$, there exists $n^*_{C;\mu,\eta}$ so that for all $n\geq n^*_{C;\mu,\eta}$
\begin{align}
N(C;\delta^*,K_n,\varepsilon^*) \geq e^{|K_n|(h_{\mu}(X,G) - \eta)}. \label{loweres1}
\end{align}

Let  $\{\xi_k\}$ be a sequence of real numbers strictly decreasing to $0$ and so that $2\xi_1 < \varepsilon^*.$ By (\ref{loweres1}), for each $k$, we find  $n_k^*$ such that for $n\geq n_k^*$ there exists a $(\delta^*,K_{n},\varepsilon^*)-$separated subset $\Gamma_{n}$ of $X_{K_{n},\mathcal{B}(\alpha_k,\xi_k/2)}$ with 
\begin{align}
|\Gamma_{n}| \geq e^{|K_{n}|h^*}. \label{sepratingcounting}
\end{align}

To get the fractal $Y_K $, we will deal with $\cup_{n}F_n$ instead of $F_n$. Thus we will use some tiling tricks to get a decomposition of $\cup_{n}F_n.$ 
We call the elements in the tempered F{\o}lner sequence $\{K_n\}$ and their translations {\itshape small bricks}.

Let $\{\tau_k\}$ be a sequence of positive numbers decreasing to $0$ satisfying $5\tau_1 < \varepsilon^*$ and $\tau_{k+1}<\tau_k/2$ for $k\geq 1.$ For each $k$, let $F(\frac{\tau_k}{2})$ be the finite subset of $G$ as described in the specification property. We assume that the identity $e_{G}$ of $G$ belongs to $F(\frac{\tau_k}{2})$ and $F(\frac{\tau_k}{2})=F(\frac{\tau_k}{2})^{-1}$.

Let $\{\gamma_k\}$ be a sequence of positive real numbers strictly decreasing to $0$ such that $\gamma_k < \min\{\frac{\delta^*}{12}, \frac{1}{12},\frac{\xi_k}{12}\}$.  We also assume $n_k^*$ increases so fast that we can find integers $n_{k}^{*} < n_{k,1} < n_{k,2}< \dots < n_{k,t_k} < n_{k+1}^*$ such that every $D$ which is $(K_{n_k,j},(\frac{\gamma_k}{|F(\frac{\tau_k}{2})|})^{t_k})$-invariant, $j=1,2,\dots,t_k$ can be $\frac{\gamma_k}{|F(\tau_{k}/2)|}-$quasi tiled by $K_{n_k,1}, K_{n_k, 2},\dots, K_{n_k,t_k}$.

We may also assume that $\{n_k^*\}$ increases so fast such that for $n\geq n_{k}^{*}$, $K_n$ is $(F(\frac{\tau_k}{2}),\frac{\gamma_k}{|F(\frac{\tau_k}{2})|})-$invariant.

Thus by Lemma \ref{newtiling}, there exist a congruent sequence of tilings $\{\mathcal{T}_k\}$(whose shape is denoted by $\{\mathcal{S}_k\}$) and the following properties hold: any $S \in \mathcal{S}_k$ can be $\frac{\gamma_k}{|F(\frac{\tau_k}{2})|}-$quasi tiled by $K_{n_k,1}, K_{n_k, 2},\dots, K_{n_k,t_k}$ with tiling centers $\{C_{k,S,1},C_{k,S,2},\dots, C_{k,S,t_k}\}$. We denote the $\frac{\gamma_k}{|F(\frac{\tau_k}{2})|}-$quasi-tiling of $S$ by	
	\[
	\mathcal{P}_{S} =\{K_{n_{k},i} c_{k,S,i} \mid 1\leq i \leq t_k,\ c_{k,S,i} \in C_{k,S,i}\}.
	\]
We call every $S\in\{\mathcal{S}_k\}$, $k\in\mathbb{N}$ and its translations  {\itshape standard bricks}. 	

Now we will modify the original bricks to get a pairwise disjoint tile of one standard brick which also satisfies the properties below. 

\begin{lemma}\label{section3modify}
	 Assume $k\in\mathbb{N}$ and $S\in\mathcal{S}_k.$ For each $c_{k,S,i}\in C_{k,S,i},$ there exists a subset $\tilde{T}_{c_{k,S,i}}\subset K_{n_k,i}$ such that if we   denote $\tilde{\mathcal{T}}_{k,S}=\{\tilde{T}_{c_{k,S,i}}c_{k,S,i} \mid c_{k,S,i}\in C_{k,S,i}, \ i=1,\dots,t_k\}$  and $\tilde{S} = \cup\tilde{\mathcal{T}}_{k,S}$ then the following properties holds:
	
	\begin{enumerate}

     	\item for $K^{\prime}\neq\ K^{\prime\prime}\in\tilde{\mathcal{T}}_{k,S}$, $F(\frac{\tau_k}{2})K^{\prime} \cap K^{\prime\prime} = \emptyset$;

		\item elements in $\tilde{\mathcal{T}}_{k,S}$ are pairwise disjoint and $|\tilde{T}_{c_{k,S,i}}| > (1-3\gamma_k)|K_{n_{k},i}|$;

		\item  $\tilde{S} \subset S, |\tilde{S}| > (1-4\gamma_k)|S|$ and $F(\frac{\tau_k}{2})\tilde{S}\subset S$;

		\item if $\Gamma_{c_{k,S,i}}$ is a $(\frac{\delta^*}{2}, \tilde{T}_{c_{k,S,i}}, \varepsilon^*)-$separated subset of $X_{\tilde{T}_{c_{k,S,i}},\mathcal{B}(\alpha_k,\xi_k)}$ with the maximal cardinality, then $|\Gamma_{c_{k,S,i}}|  \geq e^{|\tilde{T}_{c_{k,S,i}}|h^*}.$
	\end{enumerate}
\end{lemma}

\begin{proof}
	By Remark 2.1, for each $c_{k,S,i},$ we can find $T^{\prime}_{c_{k,S,i}}\subset K_{n_k,i}$ such that $|T^{\prime}_{c_{k,S,i}}|>(1-\frac{\gamma_k}{|F(\frac{\tau_{k}}{2})|})|K_{n_k,i}|$ and $\{T^{\prime}_{c_{k,S,i}}c_{k,S,i}\}$ are pairwise disjoint.
	
	 Define $\tilde{T}_{k,S,i} := \cap_{s\in F(\frac{\tau_k}{2})}s^{-1}T^{\prime}_{c_{k,S,i}}$. Then
	 \begin{align}
	 K_{n_k,i}\setminus \tilde{T}_{k,S,i} &\subset \{t\in K_{n_k,i} : \exists s\in F(\frac{\tau_k}{2}) \ \text{such that} \ st\notin T^{\prime}_{c_{k,S,j}}\} \cup (K_{n_k,i}\setminus T^{\prime}_{c_{k,S,i}}) \notag\\
	 &\subset \bigcup_{s\in F(\frac{\tau_k}{2})} \{t\in K_{n_k,i} :  st\notin T^{\prime}_{c_{k,S,j}}\} \cup (K_{n_k,i}\setminus T^{\prime}_{c_{k,S,i}}).\notag
	 \end{align}
Thus we have
\begin{align*}
|K_{n_k,i}\setminus \tilde{T}_{k,S,i}| &\leq \sum_{s\in F(\frac{\tau_k}{2})}\left(\big|\{t\in K_{n_k,i} :  st\notin K_{n_k,i}\}| + |\{t\in K_{n_k,i} :  st\in K_{n_k,i}\setminus T^{\prime}_{c_{k,S,i}}\}|\right)\\
&\quad + |(K_{n_k,i}\setminus T^{\prime}_{c_{k,S,i}})|\\
&< |F(\frac{\tau_k}{2})|(\frac{\gamma_{k}}{|F(\frac{\tau_k}{2})|}+\frac{\gamma_{k}}{|F(\frac{\tau_k}{2})|})|K_{n_k,i}| + |(K_{n_k,i}\setminus T^{\prime}_{c_{k,S,i}})|\\
 &<3\gamma_{k}|K_{n_k,i}|.
\end{align*} 
Then statements (1) and (2) hold.

To prove statement (3) :	 note that statements (1) and (2) and the fact that $S$ is $\gamma_k-$quasi tiled by $K_{n_k,1},\dots,K_{n_k,t_k}$, we have
	 \begin{align*}
	|\tilde{S}| &= \Big|\bigcup_{i=1}^{t_k}\bigcup_{c_{k,S,i}\in C_{k,S,i}}\tilde{T}_{c_{k,S,i}}c_{k,S,i}\Big| \\
	&= \sum_{i=1}^{t_k}\sum_{c_{k,S,i}\in C_{k,S,i}}|\tilde{T}_{c_{k,S,i}}|\geq (1-3\gamma_k)\sum_{i=1}^{t_k}\sum_{c_{k,S,i}\in C_{k,S,i}}|K_{n_k,i}|\\
	&\geq (1-3\gamma_k)(1-\gamma_k)|S|\geq (1-4\gamma_k)|S|.
	 \end{align*}
From the construction of $\tilde{T}_{c_{k,S,i}}$, we have $\tilde{T}_{c_{k,S,i}}c_{k,S,i}\subset K_{n_k,i}c_{k,S,i}\subset S.$
Thus $F(\frac{\tau_k}{2})\tilde{S}\subset S.$\\
We prove statement (4).  By (\ref{sepratingcounting}), there exists a $(\delta^*,K_{n_k,i},\varepsilon^*)-$separated subset $\Gamma_{n_k,i}$ of $X_{K_{n_k,i},\mathcal{B}(\alpha_k,\xi_k/2)}$ with $|\Gamma_{n_k,i}| \geq e^{|K_{n_k,i}|h^*}$.  Note that  $\gamma_k < \min\{\frac{\delta^*}{12}, \frac{1}{12},\frac{\xi_k}{12}\},\ \tilde{T}_{c_{k,S,i}}\subset K_{n_k,i}$ and $|\tilde{T}_{c_{k,S,i}}| > (1-3\gamma_k)|K_{n_{k,i}}|$, by Lemma \ref{eslemma1} we have $\Gamma_{n_k,i} \subset X_{\tilde{T}_{c_{k,S,i}},\mathcal{B}(\alpha_k,\xi_k)}$ and $\Gamma_{n_k,i}$ is a $(\frac{\delta^*}{2},\tilde{T}_{c_{k,S,i}},\varepsilon^*)-$separated set.  Let $\Gamma_{c_{k,S,i}}$ be a  $(\frac{\delta^*}{2}, \tilde{T}_{c_{k,S,i}}, \varepsilon^*)-$separated subset of $X_{\tilde{T}_{c_{k,S,i}},\mathcal{B}(\alpha_k,\xi_k)}$ with the maximal cardinality. Then $
|\Gamma_{c_{k,S,i}}| \geq |\Gamma_{n_k,i}| \geq e^{|\tilde{T}_{c_{k,S,i}}|h^*}.$
\end{proof}

Now we will use the standard bricks to build $\cup_{n=1}^{\infty}F_n.$ We will use some ideas from \cite[Section 3]{zhang2018topological}.

Let $\{\beta_k>0\}$ be a  sequence of real numbers strictly decreasing to $0$. We will choose an increasing sequence
$M(0) < M(1) < M(2) < \cdots$ of integers and a sequence $H(0) \subset H(1) \subset H(2) \subset H(3) \subset \cdots $ of sets in the following way.
\begin{enumerate}
	\item Let $H(0) = \emptyset$. Choose $M(0) >0$ such that $F_n$ is $\big(\cup\mathcal{S}_1, \frac{\beta_1}{|\cup\mathcal{S}_1|}\big)-$invariant for every $n\geq M(0)$.
	
	\item Choose $M(1) > M(0)$ such that for every $n \geq M(1)$, $F_n$ is $\big(\cup\mathcal{S}_2, \frac{\beta_2}{|\cup\mathcal{S}_2|}\big)-$invariant. Let $\tilde{F}_1 = \bigcup_{i = M(0) +1}^{M(1)} F_i$, $\tilde{\mathcal{T}_2} =\{ T\in \mathcal{T}_2\mid T\cap \tilde{F}_1 \neq \emptyset\}$
	and $H(1) = \cup \tilde{\mathcal{T}}_2$.
	
	\item Choose $M(2) > M(1)$ such that for every $n \geq M(2)$, $F_n$ is $\big(\cup\mathcal{S}_3, \frac{\beta_3}{|\cup\mathcal{S}_3|}\big)-$invariant and $|H(1)| < \beta_3 |F_n|$. Let $\tilde{F}_2 = \cup_{i = M(0) +1}^{M(2)} F_i$, $\tilde{\mathcal{T}_3} =\{ T\in \mathcal{T}_3\mid T\cap \tilde{F}_2 \neq \emptyset\}$
	and $H(2) = \cup \tilde{\mathcal{T}_3}$.
	
	\item Assume that $M(0) < M(1) < \dots < M(k-1)$ and $H(0) \subset H(1)\subset \dots \subset H(k-1)$ have been chosen, then choose $M(k) > M(k-1)$ such that for every $n \geq M(k)$, $F_n$ is 
	$\big(\cup\mathcal{S}_{k+1}, \frac{\beta_{k+1}}{|\cup\mathcal{S}_{k+1|}|}\big)-$invariant and 
	$|H(k-1)| < \beta_{k+1} |F_n|$. Let $\tilde{F}_k = \bigcup_{i = M(0) +1}^{M(k)} F_i$, 
$\tilde{\mathcal{T}}_{k+1} =\{ T\in \mathcal{T}_{k+1}\mid T\cap \tilde{F}_{k} \neq \emptyset\}$
	and $H(k) = \cup\tilde{\mathcal{T}}_{k+1}$.
			
\end{enumerate}
For $k \geq 2,$ denote
\[
H^{\prime}_k := \{T\in\mathcal{T}_k \mid T\subset H(k)\setminus H(k-1)\} \text{ and } H^{\prime}(k) := \cup H^{\prime}_k.
\]
Then $H(k)\setminus H(k-1) =H^{\prime}(k)$ for $k\geq 1.$

Next we will use standard bricks to cover each $F_n.$ The following lemma shows that most part of $F_n$ can be covered by standard bricks.
\begin{lemma}[Lemma 3.5, \cite{zhang2018topological}]\label{section3lemma2}
	For any $k$ and $M(k-1) < n\leq M(k),$ let 
	\[
	\Lambda_{n}^{1} = \big\{T\in H^{\prime}_k\mid  T \subset F_n \} \text{ and }
	\Lambda_{n}^{2} = \{T\in H^{\prime}_{k-1}\mid T\subset F_n\}.
	\]
Let $\Lambda_{n} = \Lambda_{n}^{1} \cup \Lambda_{n}^{2}$ and $F_n^{\prime} = \cup \Lambda_{n}.$ Then $
F_{n}^{\prime}\subset F_n \text{ and } |F_{n}^{\prime}| > (1-2\beta_k)|F_n|.
$
\end{lemma}

We will construct a subset of Cantor type, which will be denoted by $Y_K$.

For each $k \geq 1$ and $S\in\mathcal{S}_k$, define
\begin{align}\label{gammaS}
\Gamma(S) &:= \prod_{i=1}^{t_k}\prod_{c_{k,S,i}\in C_{k,S,i}}\Gamma_{c_{k,S,i}}\notag\\
&:= \{\vec{x}=(x_{c_{k,S,i}}) : x_{c_{k,S,i}} \in \Gamma_{c_{k,S,i}}, c_{k,S,i}\in C_{k,S,i}, i = 1,2,\dots,t_k\}.
\end{align}
For $Sd\in \mathcal{T}_k,$ define $\Gamma(Sd)=\Gamma(S).$
For $\vec{x}=(x_{c_{k,S,i}})\in \Gamma(S)$ and $r>0$, define
\begin{align}
B(S,\vec{x},r) := \{y\in X \mid \rho_{\tilde{T}_{c_{k,S,i}}c_{k,S,i}}(y,x_{c_{k,S,i}}) \leq r \text{ for all }\tilde{T}_{c_{k,S,i}}c_{k,S,i}\in\tilde{\mathcal{T}}_{k,S}\},
\end{align}
where $\tilde{\mathcal{T}}_{k,S}=\{\tilde{T}_{c_{k,S,i}}c_{k,S,i} \mid c_{k,S,i}\in C_{k,S,i}, \ i=1,\dots,t_k\}$ is as defined in Lemma \ref{section3modify}.
\noindent Define
\begin{align*}
\mathcal{F}_k &:= \bigcup_{Sd\in H_{k}^{\prime}}\{\tilde{T}_{c_{k,S,i}}c_{k,S,i}d \mid c_{k,S,i}\in C_{k,S,i}, i=1,\dots,t_k\};\\
\Gamma_k &:= \prod_{Sd\in H_{k}^{\prime}}\Gamma(Sd);\\
\mathcal{H}_k &:= \bigcup_{t=1}^{k} \mathcal{F}_t.
\end{align*}
From Lemma \ref{section3modify}, for $\tilde{T}_{c_{k+1,S,i}}c_{k+1,S,i}d\in \mathcal{F}_{k+1}$, we have
 \[F(\frac{\tau_{k+1}}{2})\tilde{T}_{c_{k+1,S,i}}c_{k+1,S,i}d \subset Sd\subset H(k+1)\setminus H(k),\]
 which implies
\begin{equation}
\big(F(\frac{\tau_{k+1}}{2})\tilde{T}_{c_{k+1,S,i}}c_{k+1,S,i}d\big) \cap H(k) = \emptyset. \label{empty5.1}
\end{equation}
Take $\vec{x}=(\vec{x}_{Sd}) \in \Gamma_k$, where $\vec{x}_{Sd}=(x_{c_{k,S,i}}) \in \Gamma(Sd)$. Choose $y\in X$ such that
\[
\rho_{\tilde{T}_{c_{k,S,i}}c_{k,S,i}d}(y,c_{k,S,i}^{-1}d^{-1}x_{c_{k,S,i}}) \leq \tau_k/2, \text{ for all }\tilde{T}_{c_{k,S,i}}c_{k,S,i}d \in \mathcal{F}_k.
\]
Such $y$ exists since the system has the specification property.

\noindent For eack $k\in\mathbb{N}$, define 
$
D_k := \{y=y(\vec{x}) \mid \vec{x} \in \Gamma_k\}.
$
Remark that 
\begin{align}
D(\mathcal{E}_{\tilde{T}_{c_{k,S,i}}c_{k,S,i}d}(y),\alpha_k) \leq \xi_k + \tau_k/2, \text{ for all }\tilde{T}_{c_{k,S,i}}c_{k,S,i}d \in \mathcal{F}_k,\  y\in D_k.\label{remark5.1}
\end{align}
	
Let $L_1 = D_1.$ Now we will define recursively $L_k$ as follows.
Suppose we have already defined the set $L_k$.  For $x\in L_k$ and $y\in D_{k+1}$, let $z=z(x,y)\in X$ be some point such that:
\begin{enumerate}
	\item for $Sd\in \mathcal{F}_{t},\ t=1,\dots,k$,
	\begin{equation}
	\rho_{\tilde{T}_{c_{t,S,i}}c_{t,S,i}d}(z(x,y),x) \leq \tau_{k+1}/2; \label{constrction1}
	\end{equation}
	\item for $Sd\in \mathcal{F}_{k+1}$,
	\begin{equation}
	\rho_{\tilde{T}_{c_{k+1,S,i}}c_{k+1,S,i}d}(z(x,y),y) \leq \tau_{k+1}/2. \label{construction2}
	\end{equation}
\end{enumerate}
By (\ref{empty5.1}), $F(\frac{\tau_{k+1}}{2})=F(\frac{\tau_{k+1}}{2})^{-1}$ and $\bigcup \mathcal{H}_k \subset H(k)$, so we know that such $z(x,y)$ exists due to the specification property. Collect all these points into the set 
\begin{equation}
L_{k+1} = \{z=z(x,y) \mid x\in L_k \text{ and } y\in D_{k+1}\}. \label{construction3}
\end{equation}
Put
\begin{equation}
Y_k = \bigcup_{x\in L_k}\{ y\in X\mid \rho_{\tilde{T}_{c_{t,S,i}}c_{t,S,i}d}(y,x) \leq \tau_k,\text{ for }\tilde{T}_{c_{t,S,i}}c_{t,S,i}d\in\mathcal{F}_{t},\ 1\leq t \leq k\}. \label{construction4}
\end{equation}
For each $y\in Y_{k+1}$, by (\ref{constrction1})--(\ref{construction4}), we know there exists $x\in L_k$ such that
\[
\rho_{\tilde{T}_{c_{t,S,i}}c_{t,S,i}d}(y,x) \leq \tau_{k+1}/2 + \tau_{k+1} < \tau_k, \text{ for }\tilde{T}_{c_{t,S,i}}c_{t,S,i}d\in\mathcal{F}_{t},\  1\leq t \leq k,
\] 
which means $y\in Y_{k}$. Thus $Y_{k+1} \subset Y_k$.

Finally, define \[Y_K = \bigcap_{k\geq 1}Y_k.\]
It follows the construction that $Y_K \neq \emptyset.$
\begin{lemma}\label{separatedes1}
	Let $k\in\mathbb{N}$ and $S\in\mathcal{S}_k.$ Then the following hold:
	\begin{enumerate}
		\item $|\Gamma(S)| \geq e^{(1-4\gamma_k)|S|h^{*}}$.
		\item Let $\vec{x} \neq \vec{y}\in \Gamma(S)$. If $x\in B(S,\vec{x},\tau_k/2)$ and $y\in B(S,\vec{y},\tau_k/2)$, then there exists $s\in \tilde{S}$ such that
		\[
		\rho(sx,sy) \geq \varepsilon^*/2.
		\]
	
	\end{enumerate}
\end{lemma}
\begin{proof}
	
(1). By statements (2-4) of Lemma \ref{section3modify}, we have
\begin{align*}
|\Gamma(S)| &= \prod_{i=1}^{t_k}\prod_{c_{k,S,i}\in C_{k,S,i}}|\Gamma_{c_{k,S,i}}|\\
&\geq e^{\sum_{i=1}^{t_k}\sum_{c_{k,S,i}\in C_{k,S,i}}|\tilde{T}_{k,S,i}|h^{*}}\\
&\geq e^{(1-4\gamma_k)|S|h^{*}}.
\end{align*}
(2). Since $\vec{x} \neq \vec{y}\in \Gamma(S)$, there exists $1\leq i \leq t_k$ and $c_{k,S,i}\in C_{k,S,i}$ such that $x_{c_{k,S,i}} \neq y_{c_{k,S,i}}\in \Gamma_{c_{k,S,i}}$. Then there exists $s\in \tilde{T}_{c_{k,S,i}}c_{k,S,i}$ such that $\rho(sx,sy) \geq \rho(sx_{c_{k,S,i}},sy_{c_{k,S,i}}) - \tau_{k} \geq \varepsilon^*/2$.
The proof is finished.
\end{proof}
Before calculating the topological entropy of $Y_K$, we state a lemma first.


	

\begin{lemma}
		 We have $Y_K\subset G_K(\mathcal{F})$.
\end{lemma}
\begin{proof}
 We define the stretched sequence $\{\alpha^\prime_n\}$ by
	\[
	\alpha_n^\prime = \alpha_k \text{ if }M(k-1)<n\leq M(k) \text{ for some } k\in\mathbb{N}.
	\]
The sequence $\{\alpha^{\prime}_m\}$  has the same limit-point set as the sequence $\{\alpha_n\}$. 

If $
	\lim\limits_{n\rightarrow \infty}D(\mathcal{E}_{F_n}(y),\alpha_n^\prime) = 0$
	then the sequences $\{\mathcal{E}_{F_n}(y)\}$ and $\{\alpha^\prime_n\}$ have the same limit-point set. Thus it is sufficient to show that for $y\in Y_{K}$ we have
$
	\lim\limits_{n\rightarrow \infty}D(\mathcal{E}_{F_n}(y),\alpha_n^\prime) = 0.
$
	
Let $k>1$ and $M(k)<n\leq M(k+1)$.	Let $\Lambda_n = \Lambda_n^{1} \cup \Lambda_n^2$ be as described in Lemma \ref{section3lemma2}. 

 For $Sd \in \Lambda_n^1,$ by (\ref{remark5.1}), (\ref{construction3}), (\ref{construction4}) and the definition of distance defined in (\ref{metriconx}), we have
\begin{align}
D(\mathcal{E}_{\tilde{T}_{c_{k+1,S,i}}c_{k+1,S,i}d}\,y,\alpha^{\prime}_{n}) = D(\mathcal{E}_{\tilde{T}_{c_{k+1,S,i}}c_{k+1,S,i}d} \,y,\alpha_{k+1}) \leq \xi_{k+1} + 2\tau_{k+1}, \label{zhlemmaes1}
\end{align}
for each $\tilde{T}_{c_{k+1,S,i}}c_{k+1,S,i} \in \tilde{\mathcal{T}}_{k+1,S}.$

For $Sd \in \Lambda_n^2,$ we have the similar estimation,
\begin{align}
D(\mathcal{E}_{\tilde{T}_{c_{k,S,i}}c_{k,S,i}d}\,y,\alpha^{\prime}_{n}) &\leq D(\mathcal{E}_{\tilde{T}_{c_{k,S,i}}c_{k,S,i}d} \,y,\alpha_{k}) + D(\alpha_k,\alpha_{k+1}) \notag\\
&\leq \xi_{k} + 2\tau_{k} + D(\alpha_k,\alpha_{k+1}), \label{zhlemmaes2}
\end{align}
for each $\tilde{T}_{c_{k,S,i}}c_{k,S,i} \in \tilde{\mathcal{T}}_{k,S}$.

By statement(3) of Lemma \ref{section3modify}, Lemma \ref{section3lemma2}, (\ref{zhlemmaes1}) and (\ref{zhlemmaes2}), we have
\begin{align}
D(\mathcal{E}_{F_n}(y),\alpha^{\prime}_n) &\leq 4\beta_{k+1} + D(\mathcal{E}_{F_n^{\prime}}(y),\alpha^{\prime}_n) \notag \\
&\leq 4\beta_{k+1}+\sum_{Sd\in \Lambda_n}\frac{|Sd|}{|F_n^{\prime}|}\big(D(\mathcal{E}_{Sd}(y),\alpha_n^{\prime})\big)\notag\\
&\leq 4\beta_{k+1}+8\gamma_k + \xi_k + 2\tau_k + D(\alpha_k,\alpha_{k+1}) \label{zhlemmaes3}.
\end{align}
From (\ref{zhlemmaes3}), we know that $\{\alpha_n^{\prime}\}$ and $\{\mathcal{E}_{F_n}(y)\}$ have the same limit-point set, which means $y\in G_K(\mathcal{F}).$
\end{proof}
\begin{proposition}\label{mainpropersition}
	We have $  h^B_{top}(Y_K,\mathcal{F}) \geq h^*.$
\end{proposition}

\begin{proof}
	Let $\varepsilon < \frac{\varepsilon^{*}}{4},$ we will show that for any $s<h^*$ we have
	 $M(Y_K,\varepsilon,s,\mathcal{F}) \geq 1.$

Since $Y_K$ is compact, we can consider a finite cover $\mathcal{C}$ of $Y_K$ whose members are sets $B_{F_{m}}(x,\varepsilon) \in \mathcal{C}$ with $B_{F_{m}}(x,\varepsilon) \cap Y_K \neq \emptyset.$ By definition, 
	\[
	M(Y_K,\varepsilon,N,s,\mathcal{F}) = \inf_{\mathcal{C}\in \mathcal{C}_{N}(Y_K,\varepsilon,\mathcal{F})}\sum_{B_{F_n}(x,\varepsilon)\in \mathcal{C}}e^{-s|F_n|}.
	\]
For each $\mathcal{C}\in \mathcal{C}_{N}(Y_K,\varepsilon,\mathcal{F})$, we define the cover $\mathcal{C}^{\prime}$ in which each ball $B_{F_n}(x,\varepsilon)$ is replaced by $B_{F_n^{\prime}}(x,\varepsilon)$.
Fix a new $\mathcal{C}\in\mathcal{C}_{N}(Y_K,\varepsilon,\mathcal{F})$. For each $B_{F_n}(x,\varepsilon)\in \mathcal{C},$ there exists $k$ such that $M(k-1)< n \leq M(k).$ Let $m$ be the largest such $k$. Define
\[
\mathcal{W}_{m} := \prod_{t=1}^{m}\prod_{Sd\in H_t^{\prime}}\Gamma(S)\ \text{and} \ \mathcal{W}_{\Lambda_n} := \prod_{Sd\in \Lambda_{n}}\Gamma(S).
\]
Fix $B_{F_n^{\prime}}(x,\varepsilon)\in\mathcal{C}^{\prime}.$ Choose $k$ such that $M(k-1)< n \leq M(k).$ Define
\[
\tilde{k} = k, \text{ if } Sd\in \Lambda_n \text{ and } Sd \in \mathcal{T}_k;\text{ otherwise set } \tilde{k} = k-1, \text{ if } Sd\in \Lambda_n \text{ and } Sd \in \mathcal{T}_{k-1}.
\]
Pick $z\in B_{F_n^{\prime}}(x,\varepsilon) \cap Y_K$. Choose $y=(y_{Sd})=((y_{c_{\tilde{k},S,i}})_{Sd})\in \mathcal{W}_{\Lambda_n}$ such that
\[
\rho_{\tilde{T}_{c_{\tilde{k},S,i}}c_{\tilde{k},S,i}}(z,y_{c_{\tilde{k},S,i}}) \leq \frac{\tau_{\tilde{k}}}{2}.
\]By  Lemma \ref{separatedes1} (2) and the choice of $\varepsilon$, such $y$ is uniquely defined. For $M(0)<j \leq M(m)$, we say the word $v\in\mathcal{W}_{\Lambda_j}$ is a {\it prefix} of $w\in \mathcal{W}_{m}$ if $w|_{\Lambda_j}=v.$ Note that each $w\in \mathcal{W}_{\Lambda_{n}}$ is the prefix of exactly $\frac{|\mathcal{W}_m|}{|\mathcal{W}_{\Lambda_{n}}|}$ words in $\mathcal{W}_m.$ Let $\overline{\mathcal{W}}_m = \bigcup_{n=M(0)+1}^{M(m)}\mathcal{W}_{\Lambda_{n}}.$ If $\mathcal{W}\subset \overline{\mathcal{W}}_m$ contains a prefix of each word in $\mathcal{W}_m,$ then
\[
\sum_{n=M(0)+1}^{M(m)}|\mathcal{W}\cap \mathcal{W}_{\Lambda_{n}}|\cdot\frac{|\mathcal{W}_m|}{|\mathcal{W}_{\Lambda_{n}}|} \geq |\mathcal{W}_m|.
\]
Thus if $\mathcal{W}$ contains a prefix of each word of $\mathcal{W}_m$,
\begin{equation}\label{fracconuting}
\sum_{n=M(0)+1}^{M(m)}|\mathcal{W}\cap \mathcal{W}_{\Lambda_{n}}|\cdot \frac{1}{|\mathcal{W}_{\Lambda_{n}}|} \geq 1.
\end{equation}
Note that since $\mathcal{C}^{\prime}$ is a cover of $Y_K$, each point of $\mathcal{W}_m$ has a prefix associated with some $B_{F_n^{\prime}}(x,\varepsilon)$ in $\mathcal{C}^{\prime}.$ Also for $M(k-1) < n\leq M(k)$ we have
\begin{align*}
|\mathcal{W}_{\Lambda_{n}}| &= \prod_{Sd\in \Lambda_{n}}|\Gamma(S)|\\
	&\geq e^{\sum_{Sd\in \Lambda_{n}}(1-4\gamma_k)h^*|S|}\\
	&\geq e^{(1-4\gamma_k)(1-2\beta_{k})|F_n|h^*}.
\end{align*}
Note that $s<h^*$. For $k$ large enough, we have
\begin{align}
s|F_n| \leq (1-4\gamma_k)(1-2\beta_{k})|F_n|h^*. \label{lastes1}
\end{align}
So when $N>M(k)$ and $\mathcal{C} \in \mathcal{C}_{N}(Y_K,\varepsilon,\{F_n\})$, by (\ref{fracconuting}) and (\ref{lastes1}) we have
\begin{align*}
\sum_{B_{F_n}(x,\varepsilon)\in\mathcal{C}}e^{-s|F_n|} \geq \sum_{B_{F_n}(x,\varepsilon)\in\mathcal{C}}e^{-(1-4\gamma_k)(1-2\beta_{k})|F_n|h^*}\geq \sum_{B_{F_n}(x,\varepsilon)\in\mathcal{C}}\frac{1}{|\mathcal{W}_{\Lambda_{n}}|} \geq 1.
\end{align*}
Hence $M(Y_K,\varepsilon,s,\mathcal{F}) \geq 1$
which means $h_{top}^{B}(Y_K,\varepsilon,\mathcal{F}) \geq s.$ Since $h_{top}^{B}(Y_K,\varepsilon,\mathcal{F})$ does not decrease as $\varepsilon$ deceases, we have $h_{top}^{B}(Y_K,\mathcal{F}) \geq s.$ 
The proof is finished.
\end{proof}

{\itshape Proof of Theorem \ref{lowerbound}:} Since $Y_K \subset G_K(\mathcal{F})$, Theorem \ref{lowerbound} is a corollary of Proposition \ref{mainpropersition}.\hfill\qed

\subsection{Table of symbols}

To have a better understanding, we make a table of some symbols and their meanings.

\begin{figure}[htbp]
	\centering
	\begin{tabular}{|c|c|}
		$\{K_n\}$ & a tempered F{\o}lner sequence, the sets $K_n$ and their translations are called small bricks\\
		$\{F_n\}$ & a general F{\o}lner sequence\\
        $\{\mathcal{T}_k\}, \{\mathcal{S}_k\}$ & congruent tiling of $G$ and its shape, \\&elements in $\mathcal{S}_k$ and their translations are called standard bricks\\
		$\Lambda_{n}, F_n^{\prime}$ & $\Lambda_{n}$ consist of standard bricks, $F_n^{\prime} = \cup \Lambda_{n}, F_n^{\prime}\subset F_n$\\
		$\beta_k$ & $F_n$ can be $(1-2\beta_k)$ covered by standard bricks, $|F_n^{\prime}| > (1-2\beta_k)|F_n|$\\
		$\gamma_k$ & $S\in\mathcal{S}_k$ can be $\gamma_k-$quasi tiled by $\{K_{n_k,1},\dots,K_{n_k,t_k}\}$ and $\gamma_k < \min\{\frac{\delta^*}{6}, \frac{1}{12},\frac{\xi_k}{6}\}$\\
		$\tilde{T}_{c_{k,S,i}}$ & small bricks $\tilde{T}_{c_{k,S,i}}\subset K_{n_k,j}$ and $\{\tilde{T}_{c_{k,S,i}}c_{k,S,i}\}$ $1-4\gamma_k$ covers $S$\\
		$\tilde{\mathcal{T}}_{k,S}$ & $\tilde{\mathcal{T}}_{k,S}=\{\tilde{T}_{c_{k,S,i}}c_{k,S,i} \mid c_{k,S,i}\in C_{k,S,i}, \ i=1,\dots,t_k\}$\\

		$\delta^*, \varepsilon^*$ & fixed small numbers with respect to  separated sets\\
		$\tau_k$ & the shadowing size in the specification property with  $5\tau_1<\varepsilon^*$ and $\tau_{k+1}<\tau_k/2$\\
		$\Gamma_{k,S,i}$ & $(\frac{\delta^*}{2}, \tilde{T}_{c_{k,S,i}},\varepsilon^*)-$separated subset of $\mathcal{B}(\alpha_k,\xi_k)$ with the maximum cardinality\\
		
	
	\end{tabular}
\end{figure}

\section{Proof of Theorem \ref{mainthm2} and Theorem \ref{mainthm3}}

{\itshape Proof of Theorem \ref{mainthm2}:}
  Assume $\hat{X}(\varphi,\mathcal{F})$ is not an empty set. Hence for $x\in \hat{X}(\varphi,\mathcal{F})$, there exist two different measures $\mu_1,\mu_2\in V(x,\mathcal{F})$ such that
  $\int_{X}\varphi \dd \mu_1 \neq \int_{X}\varphi \dd \mu_2.$

{\itshape Case 1:} $h_{\mu_1}(X,G) = h_{\mu_2}(X,G) = h_{top}(X,G).$\\
 Let $K = \{t\mu_1 + (1-t)\mu_2 \mid t\in[0,1]\}$. Since $\int_{X}\varphi \dd \mu_1 \neq \int_{X}\varphi \dd \mu_2$ we have $G_K \subset \hat{X}(\varphi,\mathcal{F}).$ Since the entropy map $\mu \mapsto h_{\mu}(X,G)$ is affine (see \cite[Theorem 8.1]{walters2000introduction} for example), we have that for all $\mu \in K$, $h_{\mu}(X,G) = h_{top}(X,G)$. By Theorem \ref{mainthm1},
\[
h^B_{top}(G_K(\mathcal{F}),\mathcal{F}) = \inf\{h_{\mu}(X,G)\mid \mu \in K\} = h_{top}(X,G),
\]
which implies $h^B_{top}(\hat{X}(\varphi,\mathcal{F})) = h_{top}(X,G).$

{\itshape Case 2:} $h_{\mu_1}(X,G) < h_{top}(X,G).$\\
Pick $0 < \eta < (h_{top}(X,G) - h_{\mu_1}(X,G))/4.$ By the variational principle (See \cite[Theorem 9.48]{kerr2016ergodic}), there is $\nu\in M(X,G)$ such that $h_{\nu}(X,G) > h_{top}(X,G) - \eta.$ Again using the fact the entropy map $\mu \mapsto h_{\mu}(X,G)$ is affine, choose $t\in (0,1)$ so small that $h_{t\mu_1 + (1-t)\nu}(X,G) > h_{top}(X,G) - 2\eta$ and  $h_{t\mu_2 + (1-t)\nu}(X,G) > h_{top}(X,G) - 2\eta.$ Consider the connected closed set $K = \{s(t\mu_1 + (1-t)\nu) + (1-s)(t\mu_2 + (1-t)\nu) \mid s\in [0,1]\}.$ Since $\int_{X}\varphi \dd \mu_1 \neq \int_{X}\varphi \dd \mu_2$, we have $t\mu_1 + (1-t)\nu \neq t\mu_2 + (1-t)\nu$ and $G_K\subset \hat{X}(\varphi,\{F_n\}).$ By the definition of $K$, for each $\mu\in K,$ $h_{\mu}(X,G) > h_{top}(X,G) - 2\eta.$
Thus by Theorem \ref{mainthm1} we have
\[
h^B_{top}(G_K(\mathcal{F}),\mathcal{F}) = \inf\{h_{\mu}(X,G)\mid \mu\in K\} > h_{top}(X,G) -2\eta.
\]
Theorem \ref{mainthm2} has been proved. \hfill\qed\\

{\itshape Proof of Theorem \ref{mainthm3}:} For statement (1.2), set $C_{\alpha} = \{\mu\in M(X,G)\mid \int_{X}\varphi\dd\mu = \alpha\}$ and $s = \sup\{h_{\mu}(X,G)\mid \mu\in C_{\alpha}\}$. \\
If $C_{\alpha} $ is a singleton set, then $G_{C_{\alpha}} = X(\varphi,\alpha,\mathcal{F})$. By Theorem \ref{mainthm1} we see that
\[
h^B_{top}(G_{C_\alpha},\mathcal{F}) = \{h_{\mu}(X,G)\mid C_{\alpha} = \{\mu\}\}.
\]
If $C_{\alpha}$ contains at least two different points, then we can pick $\mu_1\neq \mu_2\in C_{\alpha}.$\\
{\itshape Case 1:} $h_{\mu_1}(X,G) = h_{\mu_2}(X,G) = s.$\\
 Let $K = \{t\mu_1 + (1-t)\mu_2\mid t\in[0,1]\}$. Using the same argument as in the proof of Theorem \ref{mainthm2} above, we have
 \[
 h^B_{top}(G_K(\mathcal{F}),\mathcal{F}) =  \sup\{h_{\mu}(X,G)\mid \mu\in C_{\alpha}\}.
 \]
 
\noindent {\itshape Case 2:} $h_{\mu_1}(X,G) <\sup\{h_{\mu}(X,G)\mid \mu\in C_{\alpha}\}.$\\
 Pick $0 < \eta < (s - h_{\mu_1}(X,G))/4.$ Take $\nu\in C_{\alpha}$ such that $h_{\nu}(X,G) > s - \eta.$ Again using the fact the entropy map $\mu\mapsto h_{\mu}(X,G)$ is affine, choose $t\in (0,1)$ so small that $h_{t\mu_1 + (1-t)\nu}(X,G) > s - 2\eta$ and  $h_{t\mu_2 + (1-t)\nu}(X,G) > s - 2\eta.$ Consider the connected closed set $K = \{s(t\mu_1 + (1-t)\nu) + (1-s)(t\mu_2 + (1-t)\nu) \mid s\in [0,1]\}.$ Since $\mu_1 \neq \mu_2$, we have $t\mu_1 + (1-t)\nu \neq t\mu_2 + (1-t)\nu$. By the definition of $K$, for each $\mu\in K,$ we have $\mu \in C_{\alpha}$ and $h_{\mu}(X,G) > s - 2\eta.$
 Thus by Theorem \ref{mainthm1}, we have
 \[
 h^B_{top}(G_K(\mathcal{F}),\mathcal{F}) = \inf\{h_{\mu}(X,G)\mid \mu\in K\} > s
  -2\eta.
 \]
By the choice of $\eta,$ we obtain $h_{top}^B(X(\varphi,\alpha,\{F_n\})) \geq s.$\\
Obviously, $C_{\alpha}$ is a closed subset of $M(X,G).$ Let
\[
^{C_{\alpha}}G = \{x\in X\mid \{\mathcal{E}_{F_n}(x)\} \text{ has a limit point in } C_{\alpha}\}.
\]
Then by Proposition \ref{lemmaforthm1}, 
\[
h^B_{top}(^{C_{\alpha}}G,\mathcal{F}) \leq \sup\{h_{\mu}(X,G)\mid \mu\in C_{\alpha}\}.
\]
From the definition of $ X(\varphi,\alpha,\mathcal{F})$, we have $X(\varphi,\alpha,\mathcal{F}) \subset ^{C_{\alpha}}\!\!G$.
Thus we obtain
\[
h^B_{top}(X(\varphi,\alpha,\mathcal{F})) \leq \sup\{h_{\mu}(X,G)\mid \mu\in C_{\alpha}\}. 
\]
The proof is finished.\hfill\qed

{\bfseries Acknowledgments} We thank the referee for his careful reading and useful suggestions. For example the referee. For example, the referee suggests to use denote a F{\o}lner sequence using a single letter $\mathcal{F}$ which makes the statements of main results and definition of saturated set easier to write and read. And the referee reminds us  The first author would like to thank Prof. Zhihong Xia for all the help while the author was working at SUSTech.  Part of this work was done during that period. The authors would like to thank Prof. Hanfeng Li for all the help and suggestions. The first author was supported by National Natural Science Foundation of China(11801261), the second author was supported by Natural Science Foundation of China(11671093) and the third author was supported by National Natural Science Foundation of China (11871120). 

\section{Appendix}

\begin{lemma}\label{appes1}
	Let $F\in F(G), \beta>0, k\in\mathbb{N}$ and $\mathcal{T}$ a tiling of $G$ with a shape set $\mathcal{S}=\{S_1,\dots,S_l\}$. Suppose $F$ is $(\bigcup\mathcal{S},\frac{\beta}{|\bigcup\mathcal{S}|})-$invariant. Let  $\tilde{F}=\cup\{T\in\mathcal{T}\mid T\subset F\}$. Then $
	|\tilde{F}|>(1-\beta)|F|.$
\end{lemma}
\begin{proof}
	Let $I_{F}= F\setminus \tilde{F}.$ Then
	\begin{align*}
	I_{F} &=\{g\in F\mid \exists T\in \mathcal{T} \text{ such that } g\in T \text{ and }T\cap (G\setminus F) \neq \emptyset\}\\
	&\subset \bigcup\{Sd\mid  S\in\mathcal{S}, d\in G \text{ such that } d\in \partial_{S}(F)\}\\
	&\subset \bigcup\{(\cup\mathcal{S})d\mid d\in \partial_{\cup \mathcal{S}} (F)\}.
	\end{align*}
Thus $|I_F|\leq |{\cup \mathcal{S}}||\partial_{\cup \mathcal{S}} (F)|<\beta|F|$.
\end{proof}

%

The proof of the following proposition is inspired by the proofs of \cite[Theorem B]{eizenberg1994large} and \cite[Proposition 2.3]{pfister2004large}.
\begin{proposition}\label{closedinvariant}
	Suppose $(X,G)$ be a dynamical system. Let $\{K_n\}$ be a tempered F{\o}lner sequence and $\mu\in M(X,G)$. Suppose the system has the specification property and $\mu$ verifies the conclusion of Proposition \ref{sepes1}. Let $0< h^{\prime} < h_{\mu}(X,G).$ Then there exists $\varepsilon^{\prime}>0$, such that for any neighborhood $C$ of $\mu$, there exists a $G-$invariant closed subset $Y\subset X$ satisfying the following properties.

	\begin{enumerate}
		\item There exists $n_{C}^{\prime} \in \mathbb{N},$ such that $\mathcal{E}_{K_n}(y) \in C$ for all $y\in Y$ and $n \geq n_{C}^{\prime}$;
		\item There exists $n_{C}^{\prime\prime} \in \mathbb{N},$ such that there exists a subset $\Gamma_n$ of $Y$ which is $(K_n,\varepsilon^{\prime})-$separated and $|\Gamma_n| \geq e^{|K_n|h^{\prime}}$ for all $n \geq n_{C}^{\prime\prime}.$
	\end{enumerate}
In particular, $h_{top}(Y,G)\geq h^{\prime}$.
\end{proposition}

\begin{proof}
Take $h^{\prime} < h^* < h_{\mu}(X,G)$. Given $C\in\mathcal{N}(\mu)$, take a $f-$neighborhood $F^{(1)}\subset C$ of $\mu$ with fixed $\{f_j,\varepsilon_j : j=1,2,\dots,p\}.$ Denote $\varepsilon_{min} = \min\{\varepsilon_j \mid j=1,\dots,p\}.$ Let $\delta^*, \varepsilon^*$ and $n^{*}_{F^{(1/5)}}$ correspond to $h^*$ in the conclusion of Proposition \ref{sepes1}. Set $n^* = n^{*}_{F^{(1/5)}}$. 
	
Because $\{f_{j} \in C(X,\mathbb{R})\}$ are uniformly continuous on $X$, there exists $\triangle>0$ such that for every $j=1,\cdots,p$ and $x,y\in X$ we have $	\triangle < \varepsilon^*/3$ and 
\[
 \rho(x,y) < \triangle \Longrightarrow \big|f_{j}(x) - f_{j}(y)\big| <\varepsilon_j/5.
\]

Let $F(\triangle)$ be as described in the specification property with respect to $\triangle$.

Let $\{\gamma_k\}$ be a sequence of real numbers strictly decreasing to $0$ with $\gamma_1 < \min\{\frac{h^* - h^{\prime}}{3h^*},\frac{\delta^*}{6},\frac{1}{12},\frac{\varepsilon_{min}}{50}\}$. Take a sequence $\{n_k^*\}$ of integers such that $K_n$ is $(F(\triangle),\frac{\gamma_k}{|F(\triangle)|})-$invariant for $n\geq n_k^{*}$ and we can find integers $n_{k}^{*} < n_{k,1} < n_{k,2}< \dots < n_{k,t_k} < n_{k+1}^*$ such that every  $D$ which is $(K_{n_k,j},(\frac{\gamma_k}{|F(\triangle)|})^{t_k})-$invariant for each $j=1,2,\dots,t_k$  can be $\frac{\gamma_k}{|F(\triangle)|}-$quasi tiled by $K_{n_{k,1}}, K_{n_{k, 2}},\dots, K_{n_{k, t_k}}$.

Take $k^{\#}$ large such that  $n_{k^{\#}}^{*} \geq n^*$.

 By Lemma \ref{newtiling}, for  each  $k$, there exists a congruent tilings $\{\mathcal{T}_k\}$ with shape sets $\{\mathcal{S}_k\}$ such that each $S\in \mathcal{S}_k$ can be $\frac{\gamma_k}{|F(\triangle)|}-$quasi tiled by $K_{n_{k,1}}, K_{n_{k, 2}},\dots, K_{n_{k,t_k}}$ with tiling centers $\{C_{n_{k,1}},\dots,C_{n_{k,t_k}}\}$. By Lemma \ref{section3modify}, for $k\geq k^{\#}$, for any $S\in \mathcal{S}_k$ there exists a collection
\begin{equation}
\mathcal{F}_{k,S}=\{T_{c_{k,S,i}}c_{k,S,i} \mid c_{k,S,i}\in C_{k,S,i}, i=1,\dots,t_k\} \label{prop7.3es7}
\end{equation}
satisfying the following:
\begin{enumerate}
	\item for each $c_{k,S,i}$ we have $T_{c_{k,S,i}} \subset K_{n_k,i}$ and $|T_{c_{k,S,i}}|>(1-3\gamma_k)|K_{n_k,i}|$;
	\item $\left|\cup\mathcal{F}_{k,S} \right| > (1-4
	\gamma_k)|S|$;
	\item for $K\neq K^{\prime} \in \mathcal{F}_{k,S}$ we have $F(\triangle)K \cap K^{\prime} = \emptyset$;
	\item for each $c_{k,S,i}$, there exists a $(\frac{\delta^*}{2},T_{c_{k,S,i}},\varepsilon^*)-$separated subset $\Gamma_{c_{k,S,i}}\subset X_{T_{c_{k,S,i}},F^{(1/5)}}$ with 
	$|\Gamma_{c_{k,S,i}}|\geq e^{h^{*}|T_{c_{k,S,i}}|}$.
\end{enumerate}

 Fix $k\geq k^{\#} $.
Consider $Z_{F^{(1)},k}^{\#}$ defined by the requirement that $x\in Z_{F^{(1)},k}^{\#}$ if and only if for all $Sd\in \mathcal{T}_k$. Then there exists
\[
\vec{x} =(x_{c_{k,S,i}}) \in \prod_{i=1}^{t_k}\prod_{c_{k,S,i}\in C_{k,S,i}}\Gamma_{c_{k,S,i}}
\]
such that 
\begin{align}
\rho_{T_{c_{k,S,i}}}(c_{k,S,i}dx,x_{c_{k,S,i}}) \leq \triangle. \label{appens1}
\end{align}
The set $Z_{F^{(1)},k}^{\#}$ is non-empty since the specification property holds.

Let $\beta>0$ with $\beta < \frac{\varepsilon_{min}}{20}$ and $\beta< \frac{h^* - h^{\prime}}{3h^*}$. Choose $m_k$ large enough for $m\geq m_k$ implies $K_m$ is $(\cup \mathcal{S}_k,\frac{\beta}{|\cup\mathcal{S}_k|})-$invariant.

Let $m\geq m_k$. Define
\[
Y_{m,k} := \{x\in X\mid sx \in X_{K_m,F^{(4/5)}}, \forall s\in G\}.
\]
By the definition, $Y_{m,k}$ is a closed $G-$invariant subset. Next we will show that $Z_{F^{(1)},k}^{\#}\subset Y_{m,k}$. 

Take $s\in G$ and let $\Lambda_{K_{m}s} =\{T\in\mathcal{T}_k\mid T\subset K_{m}s\}$ and $\widetilde{K_{m}s}=\bigcup \Lambda_{K_{m}s}.$ By Lemma \ref{appes1} we see that
 \begin{equation}\label{apes7.1}
|\widetilde{K_{m}s}|>(1-\beta)|K_m|. 
\end{equation} 
 For each $Sd\in \Lambda_{K_{m}s}$, by (\ref{appens1}), for $x\in Z_{F^{(1)},k}^{\#}$ we have
 \begin{align*}
&\big|\sum_{t\in Sd}f_{j}(tx) - |S|\langle f_{j},\mu\rangle\big|&\\
&\leq \big|\sum_{t\in Sd}f_{j}(tx) - \sum_{T_{c_{k,S,i}}c_{k,S,i}\in\mathcal{F}_{k,S}}\sum_{s\in T_{c_{k,S,i}}c_{k,S,i}}f_{j}(sdx)\big|  &\\
&+ \big|\sum_{T_{c_{k,S,i}}c_{k,S,i}\in\mathcal{F}_{k,S}}\sum_{s\in T_{c_{k,S,i}}c_{k,S,i}}f_{j}(sdx) - \sum_{T_{c_{k,S,i}}c_{k,S,i}\in\mathcal{F}_{k,S}}\sum_{t\in T_{c_{k,S,i}}}f_{j}(tx_{c_{k,S,i}})\big| \\
&+ \big|\sum_{T_{c_{k,S,i}}c_{k,S,i}\in\mathcal{F}_{k,S}}\sum_{t\in T_{c_{k,S,i}}}f_{j}(tx_{c_{k,S,i}}) - |S|\langle f_{j},\mu\rangle\big| \\
&\leq \frac{4\gamma_{k}}{|F(\triangle)|}|S| + \frac{\varepsilon_j}{5}|S| + \big|\sum_{T_{c_{k,S,i}}c_{k,S,i}\in\mathcal{F}_{k,S}}\Big(\sum_{t\in T_{c_{k,S,i}}}f_{j}(tx_{c_{k,S,i}}) - |T_{c_{k,S,i}}|\langle f_{j},\mu\rangle\Big)\big| + \frac{4\gamma_{k}}{|F(\triangle)|}|S|\\
&\leq \frac{8\gamma_{k}}{|F(\triangle)|}|S| + \frac{\varepsilon_j}{5}|S| + \frac{\varepsilon_j}{5}|S| \leq \frac{8\gamma_{k}}{|F(\triangle)|}|S| + \frac{2\varepsilon_j}{5}|S|<\frac{3\varepsilon_j}{5}|S| .
\end{align*}
Thus we have 
\begin{align}
\big|\langle f_{j},\mathcal{E}_{Sd}(x)\rangle - \langle f_{j},\mu\rangle\big| < \frac{3\varepsilon_j}{5}. \label{importantlastes1}
\end{align}
Since $|\widetilde{K_{m}s}|>(1-\beta)|K_m|$, we have

\begin{align}
\big|\langle f_{j},\mathcal{E}_{K_{m}s}(x)\rangle  - \langle f_{j},\mathcal{E}_{\widetilde{K_{m}s}}(x)\rangle \big| <2\beta. \label{prop7.3es2}
\end{align}
By (\ref{importantlastes1}) - (\ref{prop7.3es2}), we have
\begin{align}
\big|\langle f_{j},\mathcal{E}_{K_{m}s}(x)\rangle - \langle f_{j},\mu\rangle\big| &\leq  \big|\langle f_{j},\mathcal{E}_{K_{m}s}(x)\rangle  - <f_{j},\mathcal{E}_{\widetilde{K_{m}s}}(x)>\big| + \big|\langle f_{j},\mathcal{E}_{\widetilde{K_{m}s}}(x)\rangle - \langle f_{j},\mu\rangle\big|\notag\\
&< 2\beta + \sum_{Sd\in \Lambda_{K_{m}s}}\frac{|S|}{|\widetilde{K_{m}s}|}|\langle f_{j},\mathcal{E}_{Sd}(x)\rangle - \langle f_{j},\mu\rangle| \notag\\
&< 2\beta+ \frac{3\varepsilon_j}{5} < \frac{4\varepsilon_j}{5}. \label{prop7.3es3}
\end{align}
From (\ref{prop7.3es3}), we get $Z_{F^{(1)},k}^{\#}\subset Y_{m,k}$.

Define 
\[
Y := \bigcap_{m\geq m_k}Y_{m,k}.
\]
Then $Y$ is a non-empty closed $G-$invariant subset of $X.$
\\Set $n_{C}^{\prime}=m_k.$ For  $m\geq n_{C}^{\prime},$ we have $Y\subset Y_{m,k}$, which implies that for $y\in Y, \mathcal{E}_{K_m}(y)\in F^{(4/5)} \subset C$. Then statement (1) is true.

Now we prove the statement (2) of this proposition. We set $n_{C}^{\prime\prime} = m_k$ and $\varepsilon^{\prime} = \frac{\varepsilon^*}{3}$.  Let $\Lambda_{K_{n}} =\{T\in\mathcal{T}_k\mid T\subset K_{n}\}$ and $\widetilde{K_{n}}=\cup \Lambda_{K_{n}}.$ By Lemma \ref{appes1} we have
\begin{equation}
|\widetilde{K_{n}}|>(1-\beta)|K_n|. 
\end{equation} 
For $S\in\mathcal{S}_k$, denote $\Gamma(S)=\prod_{i=1}^{q_k}\prod_{c_{k,S,i}\in C_{k,S,i}}\Gamma_{c_{k,S,i}}$. Set $\Gamma(K_n) = \prod_{Sd\in \Lambda_{K_{n}}}\Gamma(S)$. 

For each $n\geq n_{C}^{\prime\prime},$ we will consider a subset $Z_{n}^{\#} \subset Z_{F^{(1)},k}^{\#}$ with the following property: for each $\vec{x}=\{x_{Sd,c_{k,S,i}}\}\in \Gamma(K_n)$, there exists exactly one point $x\in Z_{n}^{\#}$ such that
\begin{align}
\rho_{T_{c_{k,S,j}}}(c_{k,S,j}dx,x_{Sd,c_{k,S,j}}) \leq \triangle. \label{map}
\end{align}
Define a map $\varPhi$ from $\Gamma(K_n)$ to $Z_{n}^{\#}$ such that $\varPhi(\vec{x})$ satisfies (\ref{map}) for every $x\in \Gamma_{K_n}$. For $\vec{x} \neq \vec{y} \in \Gamma(K_n)$, we have 
\[
\rho_{K_n}(\varPhi(\vec{x}),\varPhi(\vec{y})) \geq \varepsilon^{*}-2\triangle>\frac{\varepsilon^*}{3}=\varepsilon^{\prime}.
\]
Then the set $Z_{n}^{\#}$ is $(K_n,\varepsilon^{\prime})-$separated.
By the definition of $\Gamma(K_n)$ and $|\Gamma_{c_{k,S,i}}|\geq e^{h^{*}|T_{c_{k,S,i}}|}$ we obtain
\begin{align}
|\Gamma(K_n)| &= \prod_{Sd\in \Lambda_{K_{n}}}\prod_{i=1}^{q_k}\prod_{c_{k,S,i}\in C_{k,S,i}}|\Gamma_{c_{k,S,i}}|\notag\\
&\geq e^{h^{*}\sum_{Sd\in\Lambda_{K_{n}}}\sum_{i=1}^{q_k}\sum_{c_{k,S,i}\in C_{k,S,i}}|T_{c_{k,S,i}}|} \notag\\
&\geq e^{h^{*}(1-\beta)(1-4\gamma_k)|K_n|} \geq e^{h^{\prime}|K_n|}.
\end{align}
Thus the statement (2) is true which implies $h_{top}(Y,G)\geq h^{\prime}$.
\end{proof}

\begin{corollary}\label{entropydenseco}
	Under the hypothesis of Proposition \ref{closedinvariant}, the measure $\mu$ is entropy-approachable by ergodic measures.
\end{corollary}

\begin{proof}
	For any neighborhood $C\subset M(X)$ of $\mu$ and $h^{\prime} < h_{\mu}(X,G),$ let $F^{(1)}\subset C$ be an $f-$neighborhood of $\mu.$ Let $\{K_n\}$ be a tempered F{\o}lner sequence. From Proposition \ref{closedinvariant}, there exists a closed $G-$invariant subset $Y\subset X$ such that $h_{top}(Y,G) \geq h^{\prime}$ and $\mathcal{E}_{K_n}(y)\in F^{(1)}$ for $n \geq n_{C}^{\prime}$ and $y\in Y.$ Then by the variational principle,
	 there exists an ergodic measure $\nu$ with $\nu(Y) = 1$ and $h_{\nu}(Y,G) \geq h^{\prime}.$ Let $y\in Y$ be a generic point for $\nu$ with respect to $\{K_n\}.$ Since $\mathcal{E}_{K_n}(y) \rightarrow \nu$ and $\mathcal{E}_{K_n}(y) \in F^{(1)}$ for $n\geq n_{C}^{\prime}$, we have $\nu\in F^{(1)}.$
\end{proof}

\begin{lemma}
	The specification property implies the set of measures satisfying Proposition \ref{sepes1} is convex.
\end{lemma}

\begin{proof}
	Given any $f$-neighborhood $F^{(1)}$ of $\mu$, consider the corresponding $f-$neighborhoods $\hat{F}^{(1)}$ and $\tilde{F}^{(1)}$ of $\hat{\mu}$ and $\tilde{\mu}$ with the same $\{f_{j},\varepsilon_j\}$. For $h^{\prime}<h^{*}<h_{\mu}(X,G)$, select $\hat{h}^{\prime}<\hat{h}^{*}<\hat{h}_{\mu}(X,G)$ and $\tilde{h}^{\prime}<\tilde{h}^{*}<\tilde{h}_{\mu}(X,G)$
such that $h^{*} = t\hat{h}^{*} + (1-t)\tilde{h}^{*}$.	
	Let $\{K_n\}$ be a tempered F{\o}lner sequence.
	Let $\hat{\delta}^*,\hat{\varepsilon}^*$ and $\hat{n}^*_{(F(1/5))}$ correspond to $\hat{h}^*$ and let $\tilde{\delta}^*,\tilde{\varepsilon}^*$ and $\tilde{n}^*_{(F(1/5))}$ correspond to $\tilde{h}^*$ in the conclusion of Proposition \ref{sepes1}. Let $\varepsilon^*=\min\{\hat{\varepsilon}^*,\tilde{\varepsilon}^*\}, \delta^*=\min\{\hat{\delta}^*,\tilde{\delta}^*\}$ and $n^*=\max\{\hat{n}^*_{F^{1/5}},\tilde{n}^*_{F^{1/5}}\}$. Then for $n\geq n^*$, there exists a $(\delta^*,K_n,\varepsilon^*)-$separated sets $\hat{\Gamma}_n$ and $\tilde{\Gamma}_n$ of $\hat{X}_{K_n,\hat{F}^{(1/5)}}$ and $ \tilde{X}_{K_n,\tilde{F}^{(1/5)}}$ respectively with 
	\begin{align}
	|\hat{\Gamma}_n| \geq e^{\hat{h}^{*}|K_n|} \text{ and } |\tilde{\Gamma}_n| \geq e^{\tilde{h}^{*}|K_n|}. \label{appen1}
	\end{align}
	
	Let $\{\mathcal{T}_k\}$ and $\{\mathcal{S}_k\}$ be as described in Proposition \ref{closedinvariant}. For $S\in \mathcal{S}_k,$ let $\mathcal{F}_{k,S}=\{T_{c_{k,S,i}}c_{k,S,i} \mid c_{k,S,i}\in C_{k,S,i}, i=1,\dots,t_k\}$ defined as (\ref{prop7.3es7}). Let $\hat{\Gamma}_{c_{k,S,i}}$ be a $(\frac{\delta^*}{2},T_{c_{k,S,i}},\varepsilon^*)$-separated subset of $\hat{X}_{T^{\prime}_{c_{k,S,i}},F^{(2/5)}}$ with the maximal cardinality and $\tilde{\Gamma}_{c_{k,S,i}}$ be a $(\frac{\delta^*}{2},T_{c_{k,S,i}},\varepsilon^*)$-separated subset of $\tilde{X}_{T^{\prime}_{c_{k,S,i}},F^{(2/5)}}$ with the maximal cardinality.
 By (\ref{appen1}) and the arguments in the proof of Lemma \ref{section3modify}, we have $|\hat{\Gamma}_{c_{k,S,i}}| > e^{|K_n|\hat{h}^{*}}$ and $|\tilde{\Gamma}_{c_{k,S,i}}| > e^{|K_n|\tilde{h}^{*}}$.
 
Take $k$ so large that $\mathcal{F}_{k,S}$ can be divided into two parts $\mathcal{F}_{k,S}^{1}$ and $\mathcal{F}_{k,S}^{2}$ satisfying:
	\begin{enumerate}
		\item $\Big|\frac{|\mathcal{F}_{k,S}^{1}|}{|\mathcal{F}_{k,S}|} - t\Big| < 4\gamma_k$;
		\item $\Big|\frac{|\mathcal{F}_{k,S}^2|}{|\mathcal{F}_{k,S}|} - (1-t)\Big| < 4\gamma_k$.
		
	\end{enumerate}

 Define
 \begin{align*}
 \Gamma(S) &:= \Bigg(\prod_{T_{c_{k,S,i}}c_{k,S,i}\in\mathcal{F}_{k,S}^{1}}\hat{\Gamma}_{c_{k,S,i}}\Bigg)\Bigg(\prod_{T_{c_{k,S,i}}c_{k,S,i}\in\mathcal{F}_{k,S}^{2}}\tilde{\Gamma}_{c_{k,S,i}}\Bigg).
 \end{align*}
Next we consider $Z_{F^{(1)},k}^{\#}$ defined as in Proposition \ref{closedinvariant} using $\Gamma(S)$ defined above.
We finish the proof by following the ideas used to prove  Proposition \ref{closedinvariant} and Corollary \ref{entropydenseco} with straightforward modifications.
\end{proof}

\begin{corollary}\label{entropydesecor}
	If $(X,G)$ has the specification property, then the set of measures in $M(X,G)$ which are entropy-approachable by ergodic measures is closed under finite convex combinations.
\end{corollary}

\noindent{\itshape Proof of Theorem \ref{approthm}:}
	Corollary \ref{entropydesecor} shows that finite convex combinations of ergodic measures are entropy-approachable by ergodic measures. We just need to show every $\mu \in M(X,G)$ can be approximated by convex combinations of ergodic measures with entropy close to that of $\mu.$  This can be achieved by using the ergodic decomposition of $\mu.$ Let $\mu = \int_{E(X,G)}m \dd\tau(m)$ be the ergodic decomposition of $\mu.$ Take $\eta>0$. Let $\alpha = \{A_1,\dots,A_p\}$ be a finite partition of $E(X,G)$ with $diam(A_i) < \eta, \ i=1,\dots,p.$ For each $A_i\in \alpha,$  pick an ergodic measure $\mu_i\in A_i$ such that
	\[
	h_{\mu_i}(X,G) \geq \frac{1}{\tau(A_i)}\int_{A_i}h_{m}(X,G) \dd\tau(m) \text{ and } D\left(\mu_i,\frac{1}{\tau(A_i)}\int_{A_i}m\dd\tau(m)\right) < \eta.
	\]
	Denote $a_i = \tau(A_i), i=1,\dots, p.$ Then
	\[
	D(\mu,\sum_{i=1}^{p}a_i\mu_i) \leq \sum_{i=1}^{p}a_iD(\frac{1}{\tau(A_i)}\int_{A_i}m\dd\tau(m),\mu_i) \leq \eta,
	\]
	and 
	\begin{align}
	h_{\mu}(X,G) = \int_{E(X,G)}h_{m}(X,G)\dd \tau(m) \leq \sum_{i=1}^{p}a_ih_{\mu_i}(X,G). \label{last}
	\end{align}
For any neighborhood $C\in \mathcal{N}(\mu)$ and $h^{\prime} < h_{\mu}(X,G)$, choosing $\eta$ small enough for $\sum_{i=1}^{p}a_i\mu_i \in C$.  By (\ref{last}), $h^{\prime} < h_{\sum_{i=1}^{p}a_i\mu_i}(X,G)$. Now Corollary \ref{entropydesecor} implies there exists an ergodic $\nu \in C$ such that $h_{\nu}(X,G) > h^{\prime}$.

\vspace{1.3em}

\noindent {\bfseries Availability of data and materials} not applicable

\section*{Declarations}

\noindent{\bfseries Conflict of Interest:} All authors disclosed no relevant relationships.

\bibliography{saturated_sets}
\end{document}